# A binary search scheme for determining all contaminated specimens

**Vassilis G. Papanicolaou**

Department of Mathematics
National Technical University of Athens,
Zografou Campus, 157 80, Athens, GREECE

papanico@math.ntua.gr

**Abstract**

Specimens are collected from $N$ different sources. Each specimen has probability $p$ of being contaminated (e.g., in the case of an infectious disease, $p$ is the prevalence rate), independently of the other specimens. In many cases group testing is applicable, namely one can take small portions from several specimens, mix them together and test the mixture for contamination, so that if the test turns positive, then at least one of the samples in the mixture is contaminated.

In this paper we give a detailed probabilistic analysis of a binary search scheme, we propose, for determining all contaminated specimens. More precisely, we study the number $T(N)$ of tests required in order to find all the contaminated specimens, if this search scheme is applied. We derive recursive and, in some cases, explicit formulas for the expectation, the variance, and the characteristic function of $T(N)$. Also, we determine the asymptotic behavior of the moments of $T(N)$ as $N \to \infty$ and from that we obtain the limiting distribution of $T(N)$ (appropriately normalized), which turns out to be normal.

## 1 Introduction

Consider $N$ containers containing samples from $N$ different sources (e.g., clinical specimens from $N$ different patients, water samples from $N$ different lakes, soil specimens, food specimens, etc.). For each sample we assume that the probability of being contaminated (by, say, a virus, a toxic substance, etc.) is $p$ and the probability that it is not contaminated is $q := 1 - p$, independently of the other samples. All $N$ samples must undergo a screening procedure (say a molecular or antibody test in the case of a viral contamination, a radiation measurement in the case of a radioactive contamination, etc.) in order to determine exactly which are the contaminated ones.

One obvious way to identify all contaminated specimens is the *individual testing*, namely to test the specimens one by one individually (of course, this approach requires $N$ tests). In this work we analyze a group testing (or pool testing) approach, which can be called "binary search scheme" and requires a random number of tests. In particular, we will see that if $p$ is not too big, typically if $p < 0.224$, then by following this scheme, the expected number of tests required in order to determine all contaminated samples can be made strictly less than $N$. There is one requirement, though, for implementing this scheme, namely that each sample can undergo many tests (or that the sample quantity contained in each container suffices for many tests).

Group testing procedures have broad applications in areas ranging from medicine and engineering to airport security control, and they have attracted the attention of researchers for over seventy-five years (see, e.g., Dorfman's seminal paper [5], as well as [8], [12], [13], etc.). In fact, recently, group testing has received some special attention, partly due to the COVID-19 crisis (see, e.g., [1], [2], [3], [7], [10], [11], and the references therein).

The contribution of the present work is a detailed quantitative analysis of a precise group testing procedure which is quite natural and easy to implement. The scheme goes as follows: First we take samples from each of the $N$ containers, mix them together and test the mixture. If the mixture is not contaminated, we know that none of the samples is contaminated. If the mixture is contaminated, we split the containers in two groups, where the first group consists of the first $\lfloor N/2 \rfloor$ containers and the second group consists of the remaining $\lceil N/2 \rceil$ containers. Then we take samples from the containers of the first group, mix them together and test the mixture. If the mixture is not contaminated, then we know that none of the samples of the first group is contaminated. If the mixture is contaminated, we split the containers of the first group into two subgroups (in the same way we did with the $N$ containers) and continue the procedure. We also apply the same procedure to the second group (consisting of $\lceil N/2 \rceil$ containers).

Another instance where the above procedure applies is the following. Suppose we have a collection of $N$ identically looking coins (or cans, or pills, etc.), in which a relatively small percentage $p$ of them are counterfeit (defective). All genuine coins have the same weight, say $w_0$, while the weight of a counterfeit is $< w_0$. If we are interested in finding all counterfeits, we can, again, apply the binary search scheme. First we weigh all coins simultaneously. If their weight is $Nw_0$, we know there are not any counterfeits.

If, however, their weight is $< Nw_0$, then there are counterfeits and, as in the case of the contaminated samples, in order to find them we continue the search by splitting the collection into two subcollections consisting of $\lfloor N/2 \rfloor$ and $\lceil N/2 \rceil$ coins respectively.
The main quantity of interest in the present paper is the number $T(N)$ of tests required in order to determine all contaminated samples, if the above binary search scheme is applied.

Following the terminology of [1] we can characterize the aforementioned procedure as an *adaptive probabilistic group testing in the linear regime* ("linear" since the average number of contaminated specimens is $pN$, hence grows linearly with $N$). "Adaptive" refers to the fact that the way we perform a stage of the testing procedure depends on the outcome of the previous stages. "Probabilistic" refers to the fact that each specimen has probability $p$ of being contaminated, independently of the other specimens and, furthermore, that the required number of tests, $T(N)$, is a random variable. Let us be a bit more specific about the latter. Suppose $\xi_1, \xi_2, \dots$ is a sequence of independent Bernoulli random variables of a probability space $(\Omega, \mathcal{F}, \mathbb{P})$, with $\mathbb{P}\{\xi_k = 1\} = p$ and $\mathbb{P}\{\xi_k = 0\} = q$, $k = 1, 2, \dots$, so that $\{\xi_k = 1\}$ is the event that the $k$-th specimen is contaminated. Then, for each $N$, the quantity $T(N)$ of our binary search scheme is a specific function of $\xi_1, \xi_2, \dots, \xi_N$, say,

$$T(N) = T_N(\xi_1, \xi_2, \dots, \xi_N). \tag{1.1}$$

In particular, is not hard to see that

$$1 = T_N(0, 0, \dots, 0) \le T_N(\xi_1, \xi_2, \dots, \xi_N) \le T_N(1, 1, \dots, 1) = 2N - 1, \tag{1.2}$$

i.e. if none of the contents of the $N$ containers is contaminated, then $T(N) = 1$, whereas if all $N$ containers contain contaminated samples, then by easy induction on $N$ we can see that $T(N) = 2N - 1$ (observe that, for $N \ge 2$, if $\xi_1 = \cdots = \xi_N = 1$ then $T(N) = T(\lfloor N/2 \rfloor) + T(\lceil N/2 \rceil) + 1$, while, of course, $T(1) = 1$). Thus, $T(N)$ can become bigger than $N$, while the deterministic way of checking of the samples one by one requires $N$ tests.
Let us mention that, for a given $p$, the optimal testing procedure, namely the procedure which minimizes the expectation of the number of tests required in order to determine (with zero-error) all contaminated specimens in an adaptive probabilistic group testing, remains unknown [10]. However, Ungar [13] has shown that if

$$p \ge p^* := \frac{3 - \sqrt{5}}{2} \approx 0.382, \tag{1.3}$$

then the optimal procedure is individual testing. In practice, though, $p$ is usually quite small, much smaller than $p^*$.
We do not claim that our proposed procedure is optimal. However, for small values of $p$ (in fact, for $p < 0.15$) the numerical evidence we have (see Tables 1 and 2 in Subsections 2.1 and 3.1 respectively) suggests that, in comparison with the optimal results obtained by Aldridge [1]-[2] via simulation, our binary scheme is near-optimal.
Let us also mention that our "binary search scheme" is different from the "binary splitting" procedures of Sobel-Groll [12] (see, also, Aldridge [1] and Hwang [8]), where at each stage a group is splitted in two subgroups whose sizes may not be equal or even nearly equal. This feature makes these binary splitting procedures very hard to analyze theoretically, while the advantage of our proposed procedure is that, due to its simplicity it allows us to calculate certain relevant quantities (e.g., the expectation $\mathbb{E}[T(N)]$ and the variance $\mathbb{V}[T(N)]$) explicitly. We even managed to determine the limiting distribution of $T(N)$ (appropriately normalized), which turns out to be normal.
Finally, let us mention that in [3] the authors consider optimality in the case where the pool sizes are powers of 3 (except possibly for an extra factor of 4 for the first pool). In the case where $N = 100$ and $p = 0.02$ their average number of tests in order to determine the contaminated specimens is 20, thus their *average-case aspect ratio* is 0.2. Tables 1 and 2 below, indicate that our proposed strategy gives comparable, if not better, results.

In Section 2 we study the special case where $N = 2^n$. Here, the random variable $T(2^n)$ is denoted by $W_n$. We present explicit formulas for the expectation (Theorem 1) and the variance (Corollary 1) of $W_n$. In addition, in formulas (2.13) and (2.50) we give the asymptotic behavior of $\mathbb{E}[W_n]$ and $\mathbb{V}[W_n]$ respectively, as $n \to \infty$. Then, we determine the leading asymptotics of all the moments of $W_n$ (Theorem 3) and from that we conclude (Corollary 2) that, under an appropriate normalization $W_n$ converges in distribution to a normal random variable. Finally, at the end of the section, in Corollary 3 and in formula (2.83) we give a couple of Law-of-Large-Numbers-type results for $W_n$.
Section 3 studies the case of an arbitrary $N$. In Theorem 4, Theorem 5, and Corollary 5 respectively we derive recursive formulas for $\mu(N) := \mathbb{E}[T(N)]$, $g(z; N) := \mathbb{E}[z^{T(N)}]$, and $\sigma^2(N) := \mathbb{V}[T(N)]$, while Corollary 4 gives an explicit, albeit rather messy, formula for $\mathbb{E}[T(N)]$. We also demonstrate (see Remark 3) the nonexistence of the limit of $\mathbb{E}[T(N)]/N$, as $N \to \infty$, which is in contrast to the special case $N = 2^n$, where the limit exists. In Subsection 3.2 we show that the moments of $Y(N) := [T(N) - \mu(N)]/\sigma(N)$ converge to the moments of the standard normal variable $Z$ as $N \to \infty$. An immediate consequence is (Corollary 7) that $Y(N)$ converges to $Z$ in distribution.
The determination of the limiting behavior of a random quantity is a fundamental issue in many probabilistic models. Thus, our result that the normalized quantity $[T(N) - \mu(N)]/\sigma(N)$ converges in distribution to the standard normal variable $Z$, and that the convergence is good to the point that we also have convergence of all moments to the corresponding moments of $Z$, gives much more information about the testing procedure than the information we obtain from the behavior of $\mathbb{E}[T(N)]$. And since the number $N$ can be of the order of $10^2$, or even $10^3$, and the explicit results can get very messy for such big values of $N$, the asymptotic behavior of $T(N)$ enables us to obtain valuable information without much effort.
At the end of the paper we have included an appendix (Section 4) containing (a) a proposition used in the justification of formula (2.17), (b) the proof of equation (2.52), and (c) a lemma with two corollaries which are used in the proofs of some results of Sections 2 and 3.

## 2 The case $N = 2^n$

We first consider the case where the number of containers is a power of 2, namely

$$N = 2^n, \qquad \text{where} \;\; n = 0, 1, \dots . \tag{2.1}$$

As we have said in the introduction, the first step is to test a pool containing samples from all $N = 2^n$ containers. If this pool is not contaminated, then none of the contents of the $N$ containers is contaminated and we are done. If the pool is contaminated, we make two subpools, one containing samples of the first $2^{n-1}$ containers and the other containing samples of the remaining $2^{n-1}$ containers. We continue by testing the first of those subpools. If it is contaminated we split it again into two subpools of $2^{n-2}$ samples each and keep going. We also repeat the same procedure for the second subpool of the $2^{n-1}$ samples.
One important detail here is that if the first subpool of the $2^{n-1}$ samples turns out not contaminated, then we are sure that the second subpool of the remaining $2^{n-1}$ samples must be contaminated, hence we can save one test and immediately proceed with the splitting of the second subpool into two subpools of $2^{n-2}$ samples each. Likewise, suppose that at some step of the procedure a subpool of $2^{n-k}$ samples is found contaminated. Then this subpool is splitted further into two other subpools, each containing $2^{n-k-1}$ samples. If the first of these subpools is found not contaminated, then we automatically know that the other is contaminated and, consequently, we can save one test.
Let $W_n$ be the number of tests required to find all contaminated samples by following the above procedure. Then, as in (1.2) we have

$$1 \le W_n \le 2^{n+1} - 1 = 2N - 1 \tag{2.2}$$

(in the trivial case $n = 0$, i.e $N = 1$, we, of course, have $W_0 = 1$).

**Example 1.** Suppose that $N = 4$ (hence $n = 2$). In this case $W_2$ can take any value between 1 and 7, except for 2. For instance in the case where the samples of the first three containers are not contaminated, while the content of the fourth container is contaminated, we have (recall (1.1)) $W_2 = T_4(0,0,0,1) = 3$. On the other hand, in the case where the samples of the first and the third container are contaminated, while the samples of the second and fourth container are not contaminated, we have $W_2 = T_4(1,0,1,0) = 7$.

## 2.1 The expectation of $W_n$

**Theorem 1.** Let $N = 2^n$ and $W_n$ as above. Then

$$\mathbb{E}[W_n] = 2^{n+1} - 1 - 2^n \sum_{k=1}^{n} \frac{q^{2^k} + q^{2^{k-1}}}{2^k}, \qquad n = 0, 1, \dots \tag{2.3}$$

(in the trivial case $n = 0$ the sum is empty, i.e. 0), where, as it is mentioned in the introduction, $q$ is the probability that a sample is not contaminated.

*Proof.* Assume $n \geq 1$ and let $D_n$ be the event that none of the $2^n$ samples is contaminated. Then

$$\mathbb{E}[W_n] = \mathbb{E}[W_n | D_n]\, \mathbb{P}(D_n) + \mathbb{E}[W_n | D_n^c]\, \mathbb{P}(D_n^c),$$

and hence

$$\mathbb{E}[W_n] = q^{2^n} + u_n \left(1 - q^{2^n}\right), \tag{2.4}$$

where for typographical convenience we have set

$$u_n := \mathbb{E}[W_n | D_n^c]. \tag{2.5}$$

In order to find a recursive formula for $u_n$ let us first consider the event $A_n$ that in the group of the $2^n$ containers none of the first $2^{n-1}$ contain contaminated samples. Clearly, $D_n \subset A_n$ and

$$\mathbb{P}(A_n | D_n^c) = \frac{\mathbb{P}(A_n D_n^c)}{\mathbb{P}(D_n^c)} = \frac{\mathbb{P}(A_n) - \mathbb{P}(D_n)}{\mathbb{P}(D_n^c)} = \frac{q^{2^{n-1}} - q^{2^n}}{1 - q^{2^n}} = \frac{q^{2^{n-1}}}{1 + q^{2^{n-1}}}. \tag{2.6}$$

Likewise, if $B_n$ is the event that in the group of the $2^n$ containers none of the last $2^{n-1}$ contains contaminated samples, then

$$\mathbb{P}(B_n | D_n^c) = \frac{q^{2^{n-1}}}{1 + q^{2^{n-1}}}. \tag{2.7}$$

Let us also notice that $A_n B_n = D_n$, hence $\mathbb{P}(A_n B_n | D_n^c) = 0$.

Now, (i) given $A_n$ and $D_n^c$ we have that $W_n \stackrel{d}{=} 1 + W_{n-1}$, where the notation $X \stackrel{d}{=} Y$ signifies that the random variables $X$ and $Y$ have the same distribution. To better justify this distributional equality, let us observe that, in the case $N = 2^n$, formula (1.1) takes the form $W_n = T(N) =$ $T_N(\xi_1, \ldots, \xi_N)$ and $W_{n-1} = T(N/2) = T_{N/2}(\xi_1, \ldots, \xi_{N/2})$. Thus, $W_{n-1} \overset{d}{=} T_{N/2}(\xi_{(N/2)+1}, \ldots, \xi_N)$. Now, $D_n = \{\xi_1 = \xi_2 = \cdots = \xi_N = 0\}$ and $A_n = \{\xi_1 = \xi_2 = \cdots = \xi_{N/2} = 0\}$. Hence, the given event $A_n D_n^c$ means that only some of the variables $\xi_{(N/2)+1}, \ldots, \xi_N$ are different from 0. Consequently, given $A_n D_n^c$ we have that $W_n \overset{d}{=} 1 + W_{n-1}$.
(ii) Similarly, given $B_n$ and $D_n^c$ we have that $W_n \overset{d}{=} 2 + W_{n-1}$.
Finally, (iii) given $(A_n \cup B_n)^c$ and $D_n^c$ we have that $W_n = 1 + W_{n-1} + \tilde{W}_{n-1}$, where, $\tilde{W}_{n-1}$ is an independent copy of $W_{n-1}$. Thus (given $D_n^c$ ) by conditioning on the events $A_n$, $B_n$ and $(A_n \cup B_n)^c$ (recall that $\mathbb{P}(A_n B_n | D_n^c) = 0$) we get, in view of (2.5), (2.6), and (2.7),

$$u_n = (1 + u_{n-1}) \frac{q^{2^{n-1}}}{1 + q^{2^{n-1}}} + (2 + u_{n-1}) \frac{q^{2^{n-1}}}{1 + q^{2^{n-1}}} + (1 + 2u_{n-1}) \frac{1 - q^{2^{n-1}}}{1 + q^{2^{n-1}}}$$

or

$$u_n = \frac{2}{1 + q^{2^{n-1}}} u_{n-1} + \frac{1 + 2q^{2^{n-1}}}{1 + q^{2^{n-1}}}. \tag{2.8}$$

Also, by replacing $n$ by $n - 1$ in (2.4) we get

$$\mathbb{E}[W_{n-1}] = q^{2^{n-1}} + u_{n-1} \left(1 - q^{2^{n-1}}\right). \tag{2.9}$$

Therefore, by using (2.8) in (2.4) and (2.9) we can eliminate $u_n$, $u_{n-1}$ and obtain

$$\mathbb{E}[W_n] = 2\,\mathbb{E}[W_{n-1}] - q^{2^n} - q^{2^{n-1}} + 1. \tag{2.10}$$

Finally, (2.3) follows easily from (2.10) and the fact that $\mathbb{E}[W_0] = 1$. ■

In view of (2.48) (given below) formula (2.3) can be also written as

$$\mathbb{E}[W_n] = 2^{n+1} - 1 - 2^{n-1} q - q^{2^n} - 3 \cdot 2^{n-1} \sum_{k=1}^{n-1} \frac{q^{2^k}}{2^k}, \qquad n \geq 0. \tag{2.11}$$

As a small illustration let us mention that in the cases $n = 1$, $n = 2$, and $n = 10$ formula (2.3) (or, equivalently, (2.11)) becomes

$$\mathbb{E}[W_1] = 3 - q - q^2, \qquad \mathbb{E}[W_2] = 7 - 2q - 3q^2 - q^4,$$

and (2.12)

$$\mathbb{E}[W_{10}] = 2047 - 512q - 768q^2 - 384q^4 - 192q^8 - 96q^{16} - 48q^{32} - 24q^{64} - 12q^{128} - 6q^{256} - 3q^{512} - q^{1024}$$

respectively. In the extreme case $q = 0$ we had seen that $W_n = 2N - 1$, hence $\mathbb{E}[W_n] = 2N - 1 = 2^{n+1} - 1$, while in the other extreme case $q = 1$ it is clear that $W_n = 1$, hence $\mathbb{E}[W_n] = 1$. Both of these values agree with the ones given by (2.3).

Let us also notice that formula (2.3) implies that for any given $q \in [0,1)$ we have

$$\mathbb{E}[W_n] = 2^n \alpha_1(q) - 1 + O\left(q^{2^n}\right), \qquad n \to \infty, \tag{2.13}$$

where

$$\alpha_1(q) := 2 - \sum_{k=1}^{\infty} \frac{q^{2^k} + q^{2^{k-1}}}{2^k} = 2 - \frac{q}{2} - \frac{3}{2}\sum_{k=1}^{\infty} \frac{q^{2^k}}{2^k}. \tag{2.14}$$

It is clear that $\alpha_1(q)$ is a power series about $q = 0$ with radius of convergence equal to 1, which is strictly decreasing on $[0,1]$, with $\alpha_1(0) = 2$ and $\alpha_1(1) = 0$ (actually, since in the power series appear only powers $q^{2^k}$, the so-called gap theorem [9] implies that the unit circle is the natural boundary of the series). Furthermore, it is not hard to see that $\alpha_1(q)$ satisfies the functional equation

$$\alpha_1\left(q^2\right) = 2\alpha_1(q) + q^2 + q - 2. \tag{2.15}$$

Closing this brief analysis of $\alpha_1(q)$, let us determine the asymptotic behavior of its derivative as $q \to 1^-$. From the second equality of (2.14) we have

$$\alpha_1'(q) = -\frac{1}{2} - \frac{3}{2q}\sum_{k=1}^{\infty} q^{2^k}, \tag{2.16}$$

thus, by using in (2.16) the inequality

$$\int_1^{\infty} q^{2^x}\, dx < \sum_{k=1}^{\infty} q^{2^k} < 1 + \int_1^{\infty} q^{2^x}\, dx, \qquad q \in (0,1),$$

together with the asymptotic estimate of the above integral given in Proposition A1 of the Appendix, we get that

$$\alpha_1'(q) = \frac{3}{2\ln 2}\ln(1-q) + O(1), \qquad q \to 1^-. \tag{2.17}$$

**Remark 1.** By dividing both sides of (2.3) by $2^n$ and comparing with (2.14) we can see that as $n \to \infty$

$$\frac{\mathbb{E}[W_n]}{2^n} \to \alpha_1(q) \qquad \text{uniformly in } \; q \in [0,1]. \tag{2.18}$$

Actually, for every $m = 1, 2, \ldots$ the $m$-th derivative of $\mathbb{E}[W_n]/2^n$ with respect to $q$ converges to the $m$-th derivative of $\alpha_1(q)$ uniformly on compact subets of $[0,1)$.

The expectation $\mathbb{E}[W_n]$ and, consequently, the *average-case aspect ratio* $\mathbb{E}[W_n]/2^n$ are getting smaller as $q$ approaches 1 (equivalently, as $p$ approaches 0). As an illustration, by employing (2.3) we have created the following table:

**Table 1:** *Values of* $\mathbb{E}[W_n]/N$ *for certain choices of* $N = 2^n$ *and* $p$.

| | $p = .005$ | $p = .01$ | $p = .05$ | $p = .1$ | $p = .15$ | $p = .2$ |
|---|---|---|---|---|---|---|
| $N = 1$ | 1 | 1 | 1 | 1 | 1 | 1 |
| $N = 2$ | .505 | .515 | .575 | .645 | .715 | .780 |
| $N = 4$ | .265 | .280 | .395 | .528 | .653 | .768 |
| $N = 8$ | .148 | .169 | .335 | .518 | .679 | .820 |
| $N = 16$ | .092 | .121 | .328 | .542 | .719 | .871 |
| $N = 32$ | .068 | .103 | .340 | .566 | .748 | .901 |
| $N = 64$ | .059 | .099 | .352 | .581 | .763 | .917 |
| $N = 128$ | .057 | .101 | .359 | .589 | .771 | .924 |
| $N = 256$ | .058 | .103 | .363 | .593 | .775 | .928 |
| $N = 512$ | .059 | .105 | .365 | .595 | .777 | .930 |
| $N = 1024$ | .060 | .106 | .366 | .596 | .778 | .931 |

(by comparing the above values with the graphs, given in [1] and [2], of the optimal average-case aspect ratio as a function of $p$, we can see that, if $p \le 0.15$, our procedure is near-optimal).

If $q$ is such that

$$\mathbb{E}[W_n] < 2^n = N, \tag{2.19}$$

then, by applying the aforementioned testing strategy, the number of tests required to determine all containers with contaminated content is, on the average, less than $N$, namely less than the number of tests required to check the containers one by one.
Let us set

$$\mu_n = \mu_n(q) := \mathbb{E}[W_n]. \tag{2.20}$$

It is clear from (2.3) that for $n \ge 1$ the quantity $\mu_n(q)$ is a polynomial in $q$ of degree $2^n$, strictly decreasing on $[0, 1]$, with $\mu_n(0) = 2^{n+1} - 1 = 2N - 1$ and $\mu_n(1) = 1$. Thus, there is a unique $q_n \in (0, 1)$ such that

$$\mu_n(q_n) = 2^n \tag{2.21}$$

and (2.19) holds if $q \in (q_n, 1]$ or, equivalently, if $p \in [0, 1 - q_n)$. Therefore, it makes sense to apply our binary search procedure only if $q > q_n$.
As an example, let us give the (approximate) values of $q_n$ for $n = 0, 1, \dots, 10$:

$$\begin{array}{lllll} q_0 = 0, & q_1 = \frac{\sqrt{5}-1}{2} \approx .618, & q_2 \approx .685, & q_3 \approx .727, & q_4 \approx .751, \\ q_5 \approx .763, \quad q_6 \approx .769, & q_7 \approx .772, & q_8 \approx .774, & q_9 \approx .775, & q_{10} \approx .775, \end{array}$$

where $q_0 = 0$ is a convention. Notice that $q_1 = 1 - p^*$, where $p^*$ is as in (1.3). Thus, $q_1$ is the threshold found by Ungar [13] and, consequently, if $q \le q_1$ the group testing does not improve on individual testing.

In view of (2.3) and (2.20), formula (2.21) becomes

$$\sum_{k=1}^{n} \frac{q_n^{2^k} + q_n^{2^{k-1}}}{2^k} = 1 - \frac{1}{2^n}, \tag{2.22}$$

which by letting $n \to \infty$ yields

$$q_n \to q_\infty, \qquad \text{where } \alpha_1(q_\infty) = 1, \tag{2.23}$$

$\alpha_1(q)$ being the function defined in (2.14). In fact, $.775 < q_\infty < .776$. Furthermore, by formula (2.22) we have

$$\sum_{k=1}^{n} \frac{q_n^{2^k} + q_n^{2^{k-1}}}{2^k} = \sum_{k=1}^{n} \frac{q_{n+1}^{2^k} + q_{n+1}^{2^{k-1}}}{2^k} - \frac{1 - q_{n+1}^{2^n} - q_{n+1}^{2^{n+1}}}{2^{n+1}},$$

which implies that the sequence $q_n$ is strictly increasing. Therefore, if $q > .776$, equivalently if $p < .224$, then $\mathbb{E}[W_n] < 2^n = N$ for any $n \geq 1$.
Finally, let us mention that, for a given $q$, it is worth to estimate the optimal value, say $m = m(q)$, of $N$ which minimizes $\mathbb{E}[T(N)]/N$, since, having $m(q)$, it may make sense to split our original $N$ samples into subgroups, each having size $m(q)$. For instance, Table 1 suggests that $m(.995) \approx 128$, while $m(.99) \approx 64$; in the case, though, where $q$ falls below .85 (see also Table 2 in Subsection 3.1), the value of $m(q)$ becomes $\leq 5$.
Actually, there is a heuristic way to estimate $m(q)$, at least in the case where $N = 2^n$: From formula (2.3) we have

$$\frac{\mathbb{E}[W_n]}{2^n} = 2 - 2^{-n} - \sum_{k=1}^{n} \frac{q^{2^k} + q^{2^{k-1}}}{2^k}, \qquad n \geq 0. \tag{2.24}$$

For large $n$ we can approximate the above expression for $\mathbb{E}[W_n]/2^n$ by the smooth function

$$f(x) = 2 - 2^{-x} - \int_0^x \frac{q^{2^\xi} + q^{2^{\xi-1}}}{2^\xi}\, d\xi, \qquad x \geq 0, \tag{2.25}$$

so that $f(n) \approx \mathbb{E}[W_n]/2^n$. Thus, in order to get an estimate of the value of $N = 2^n$ which minimizes $\mathbb{E}[W_n]/2^n$, we can consider the minimization of $f(x)$ on $[0, \infty)$. We have

$$f'(x) = 2^{-x} \ln 2 - \frac{q^{2^x} + q^{2^{x-1}}}{2^x}, \tag{2.26}$$

hence, the minimum of $f(x)$ is attained at the value of $x$ for which $f'(x) = 0$, i.e.

$$\ln 2 = q^M + q^{M/2}, \qquad \text{where } M = 2^x, \tag{2.27}$$

which implies (since $q^{M/2} > 0$)

$$q^{M/2} = \frac{\sqrt{1 + 4\ln 2} - 1}{2} \approx 0.47 \tag{2.28}$$

or

$$M = \frac{2\left[\ln\left(\sqrt{1 + 4\ln 2} - 1\right) - \ln 2\right]}{\ln q} \approx \frac{1.5}{-\ln q} \tag{2.29}$$

The interesting case is when the prevalence rate $p = 1 - q$ is quite small. In this case $-\ln q \approx p$, hence $M \approx 1.5/p$, which suggests that the $N = 2^n$ which minimizes $\mathbb{E}[T(N)]/N$ satisfies $N \sim c/p$ as $p \to 0^+$, for some constant $c > 0$.

## 2.2 The behavior of $W_n$ as $n \to \infty$

We start with a recursive formula for the generating function of $W_n$.

**Theorem 2.** Let

$$g_n(z) := \mathbb{E}\left[z^{W_n}\right]. \tag{2.30}$$

Then

$$g_n(z) = z g_{n-1}(z)^2 + (z - z^2) q^{2^{n-1}} g_{n-1}(z) + (z - z^2) q^{2^n}, \qquad n = 1, 2, \dots \tag{2.31}$$

(clearly, $g_0(z) = z$).

*Proof.* With $D_n$ as in the proof of Theorem 1 we have

$$g_n(z) = \mathbb{E}\left[z^{W_n}\right] = \mathbb{E}\left[z^{W_n} | D_n\right] \mathbb{P}(D_n) + \mathbb{E}\left[z^{W_n} | D_n^c\right] \mathbb{P}(D_n^c),$$

thus

$$g_n(z) = z q^{2^n} + h_n(z)\left(1 - q^{2^n}\right), \tag{2.32}$$

where

$$h_n(z) := \mathbb{E}\left[z^{W_n} | D_n\right]. \tag{2.33}$$

Let $A_n$ and $B_n$ be the events introduced in the proof of Theorem 1. Then, given $A_n$ and $D_n^c$ we have that $z^{W_n} \stackrel{d}{=} z z^{W_{n-1}}$; given $B_n$ and $D_n^c$ we have that $z^{W_n} \stackrel{d}{=} z^2 z^{W_{n-1}}$; finally given $(A_n \cup B_n)^c$ and $D_n^c$ we have that $z^{W_n} \stackrel{d}{=} z z^{W_{n-1}} z^{\tilde{W}_{n-1}}$, where $\tilde{W}_{n-1}$ is an independent copy of $W_{n-1}$. Thus (given $D_n^c$ ) by conditioning on the events $A_n$, $B_n$ and $(A_n \cup B_n)^c$ we get, in view of (2.33), (2.6), and (2.7),

$$h_n(z) = (z + z^2) h_{n-1}(z) \frac{q^{2^{n-1}}}{1 + q^{2^{n-1}}} + z h_{n-1}(z)^2 \frac{1 - q^{2^{n-1}}}{1 + q^{2^{n-1}}}. \tag{2.34}$$

By replacing $n$ by $n - 1$ in (2.32) we get

$$g_{n-1}(z) = z q^{2^{n-1}} + h_{n-1}(z)\left(1 - q^{2^{n-1}}\right) \tag{2.35}$$

and (2.31) follows by eliminating $h_n(z)$ and $h_{n-1}(z)$ from (2.32), (2.34), and (2.35). ■

By setting $z = e^{it}$ in (2.31) it follows that the characteristic function

$$\phi_n(t) := \mathbb{E}\left[e^{itW_n}\right] \tag{2.36}$$

of $W_n$ satisfies the recursion

$$\phi_n(t) = e^{it}\phi_{n-1}(t)^2 + \left(e^{it} - e^{2it}\right) q^{2^{n-1}} \phi_{n-1}(t) + \left(e^{it} - e^{2it}\right) q^{2^n}, \quad n = 1, 2, \dots, \tag{2.37}$$

with, of course, $\phi_0(t) = e^{it}$.

For instance, for $n = 1$ formula (2.37) yields

$$\phi_1(t) = e^{3it} + q\left(e^{it} - e^{2it}\right) e^{it} + q^2\left(e^{it} - e^{2it}\right) e^{it}. \tag{2.38}$$

**Corollary 1.** For the variance

$$\sigma_n^2 = \sigma_n^2(q) := \mathbb{V}\left[W_n\right] \tag{2.39}$$

of $W_n$ we have the formula

$$\sigma_n^2(q) = 2^{n+1}\sum_{k=1}^{n}\left(2q^{2^k} + q^{2^{k-1}}\right) + 2^n \sum_{k=1}^{n} \frac{q^{2^{k+1}} + q^{3\cdot 2^{k-1}} - 5q^{2^k} - 3q^{2^{k-1}}}{2^k}$$
$$- 2^n \sum_{k=1}^{n}\left(2q^{2^k} + q^{2^{k-1}}\right)\sum_{j=1}^{k}\frac{q^{2^j} + q^{2^{j-1}}}{2^j}, \qquad n = 0, 1, \dots \tag{2.40}$$

(in the trivial case $n = 0$ all the above sums are empty, i.e. 0).

*Proof.* Differentiating (2.31) twice with respect to $z$ and then setting $z = 1$ yields

$$g_n''(1) = 2g_{n-1}''(1) + 2g_{n-1}'(1)^2 - 2q^{2^{n-1}} g_{n-1}'(1) + 4g_{n-1}'(1) - 2q^{2^n} - 2q^{2^{n-1}}, \tag{2.41}$$

where, in view of (2.30),

$$g_n''(1) = \mathbb{E}\left[W_n(W_n - 1)\right] \qquad \text{and} \qquad g_n'(1) = \mathbb{E}\left[W_n\right]. \tag{2.42}$$

Thus,

$$\sigma_n^2 = g_n''(1) + g_n'(1) - g_n'(1)^2 \tag{2.43}$$

and by using (2.10) and (2.43) in (2.41) we obtain

$$\sigma_n^2 = 2\sigma_{n-1}^2 + \left(2q^{2^n} + q^{2^{n-1}}\right)\mathbb{E}\left[W_n\right] + q^{2^{n+1}} + q^{3\cdot 2^{n-1}} - 3q^{2^n} - 2q^{2^{n-1}}. \tag{2.44}$$

Finally, formula (2.44) together with the fact that $\sigma_0^2 = \mathbb{V}\left[W_0\right] = 0$ imply

$$\sigma_n^2 = 2^n \sum_{k=1}^{n} \frac{\left(2q^{2^k} + q^{2^{k-1}}\right)\mathbb{E}[W_k] + q^{2^{k+1}} + q^{3\cdot 2^{k-1}} - 3q^{2^k} - 2q^{2^{k-1}}}{2^k} \tag{2.45}$$

from which (2.40) follows by invoking (2.3). ■

As we have mentioned, in the extreme cases where $q = 0$ or $q = 1$ the variable $W_n$ becomes deterministic and, consequently, $\sigma_n^2(0) = \sigma_n^2(1) = 0$, which is in agreement with (2.40).

In the cases $n = 1$ and $n = 2$ formula (2.40) yields respectively

$$\sigma_1^2(q) = \mathbb{V}[W_1] = -q^4 - 2q^3 + 2q^2 + q = q(1-q)(q^2 + 3q + 1) \tag{2.46}$$

and

$$\sigma_2^2(q) = \mathbb{V}[W_2] = q(1-q)(q^6 + q^5 + 7q^4 + 11q^3 + 5q^2 + 11q + 2). \tag{2.47}$$

Notice that the maximum of $\sigma_1^2(q)$ on $[0,1]$ is attained at $q \approx .6462$. Actually, both $d\sigma_1^2(q)/dq$ and $d^2\sigma_1^2(q)/dq^2$ have one (simple) zero in $[0,1]$. The maximum of $\sigma_2^2(q)$ on $[0,1]$ is attained at $q \approx .744$, while both $d\sigma_2^2(q)/dq$ and $d^2\sigma_2^2(q)/dq^2$ have one simple zero in $[0,1]$.

Let us also notice that, since

$$\sum_{j=1}^{k} \frac{q^{2^j} + q^{2^{j-1}}}{2^j} = \frac{q}{2} + \frac{q^{2^k}}{2^k} + \frac{3}{2}\sum_{j=1}^{k-1} \frac{q^{2^j}}{2^j}, \tag{2.48}$$

formula (2.40) can be also written as

$$\sigma_n^2(q) = 2^n \left(2 - \frac{q}{2}\right) \sum_{k=1}^{n} \left(2q^{2^k} + q^{2^{k-1}}\right) - 2^n \sum_{k=1}^{n} \frac{q^{2^{k+1}} + 5q^{2^k} + 3q^{2^{k-1}}}{2^k}$$
$$- 2^n \cdot \frac{3}{2} \sum_{k=1}^{n} \left(2q^{2^k} + q^{2^{k-1}}\right) \sum_{j=1}^{k-1} \frac{q^{2^j}}{2^j}, \qquad n = 0, 1 \ldots . \tag{2.49}$$

From formula (2.40) it follows that for any given $q \in [0,1)$ we have

$$\sigma_n^2 = \mathbb{V}[W_n] = 2^n \beta(q) + O\left(2^n q^{2^n}\right), \qquad n \to \infty, \tag{2.50}$$

where

$$\beta(q) := \sum_{k=1}^{\infty} \left(2 - \sum_{j=1}^{k} \frac{q^{2^j} + q^{2^{j-1}}}{2^j}\right) \left(2q^{2^k} + q^{2^{k-1}}\right)$$
$$+ \sum_{k=1}^{\infty} \frac{q^{2^{k+1}} + q^{3 \cdot 2^{k-1}} - 5q^{2^k} - 3q^{2^{k-1}}}{2^k}. \tag{2.51}$$

Two equivalent expressions for $\beta(q)$ are

$$\beta(q) = \left(2 - \frac{q}{2}\right) \sum_{k=1}^{\infty} \left(2q^{2^k} + q^{2^{k-1}}\right) - \sum_{k=1}^{\infty} \frac{q^{2^{k+1}} + 5q^{2^k} + 3q^{2^{k-1}}}{2^k}$$
$$- \frac{3}{2} \sum_{k=1}^{\infty} \left(2q^{2^k} + q^{2^{k-1}}\right) \sum_{j=1}^{k-1} \frac{q^{2^j}}{2^j}$$

and

$$\beta(q) = \frac{q(q+1)}{2} + \left(6 - \frac{3q}{2}\right) \sum_{k=1}^{\infty} q^{2^k} - \frac{17}{2} \sum_{k=1}^{\infty} \frac{q^{2^k}}{2^k}$$
$$- \frac{3}{2} \sum_{k=1}^{\infty} \left(2q^{2^k} + q^{2^{k-1}}\right) \sum_{j=1}^{k-1} \frac{q^{2^j}}{2^j}$$

(in the power series of $\beta(q)$ appear only powers $q^{n_k}$, where $n_k = 2^k$, or $n_k = 2^k + 1$, or $n_k = 2^k + 2$, or $n_k = 2^k + 4$, $\dots$, or $n_k = 2^k + 2^{k-1}$, hence the so-called gap theorem [9] implies that the unit circle is the natural boundary of the series).

From (2.51) we can see that $\beta(q)$ is finite for $q \in [0, 1]$, with $\beta(0) = \beta(1) = 0$. In fact we can say a bit more, namely

$$\lim_{q \to 1^-} \beta(q) = \beta(1) = 0 \tag{2.52}$$

(the proof of (2.52) is given in the Appendix). Let us also notice that from formula (2.50) we can see that $\beta(q) > 0$ for $q \in (0, 1)$ since it is related to the variance of $W_n$.

We would like now to understand the limiting behavior of $W_n$, as $n \to \infty$, for any fixed $q \in (0, 1)$. Formulas (2.13) and (2.50) suggest that we should work with the normalized variables

$$X_n := \frac{W_n - 2^n \alpha_1(q)}{2^{n/2}}, \qquad n = 0, 1, \dots . \tag{2.53}$$

Let us consider the characteristic functions

$$\psi_n(t) := \mathbb{E}\left[e^{itX_n}\right], \qquad n = 0, 1, \dots . \tag{2.54}$$

In view of (2.36) we have

$$\psi_n(t) = \phi_n\left(2^{-n/2} t\right) e^{-i 2^{n/2} \alpha_1(q) t} \tag{2.55}$$

and using (2.55) in (2.37) yields

$$\begin{aligned} \psi_n(t) &= e^{2^{-n/2}\, it}\, \psi_{n-1}\left(t/\sqrt{2}\right)^2 \\ &+ \left(e^{2^{-n/2}\, it} - e^{2 \cdot 2^{-n/2}\, it}\right) e^{-2^{n/2}\, \alpha_1(q) it/2} q^{2^{n-1}} \psi_{n-1}\left(t/\sqrt{2}\right) \\ &+ \left(e^{2^{-n/2}\, it} - e^{2 \cdot 2^{-n/2}\, it}\right) e^{-2^{n/2}\, \alpha_1(q) it} q^{2^n}, \qquad n = 1, 2, \dots, \end{aligned} \tag{2.56}$$

with $\psi_0(t) = e^{[1 - \alpha_1(q)] it}$.

From formula (2.56) one can guess the limiting distribution of $X_n$: Assuming that for each $t \in \mathbb{R}$ the sequence $\psi_n(t)$ has a limit, say $\psi(t)$, we can let $n \to \infty$ in (2.56) and conclude that $\psi(t) = \psi\left(t/\sqrt{2}\right)^2$. Actually, if the limit $\psi(t)$ exists, then formulas (2.57) and (2.58) below imply that $\psi'(0) = 0$ and $\psi''(0)$ exists (in $\mathbb{C}$). From the existence of these derivatives and the functional equation of $\psi(t)$ we get that the function $h(t) := t^{-1}\sqrt{-\ln\psi(t)}$ is continuous at 0 and satisfies the equation $h(t) = h\left(t/\sqrt{2}\right)$. Therefore, $h(t) \equiv h(0)$ for all $t \in \mathbb{R}$ and, hence $\psi(t) = e^{-at^2}$ for some $a \geq 0$. Consequently, $X_n$ converges in distribution to a normal random variable, say $X$, with zero mean. As for the variance of $X$, formulas (2.50) and (2.53) suggest that $\mathbb{V}[X] = \beta(q)$, where $\beta(q)$ is given by (2.51).

The above argument is not rigorous since it assumes the existence of the (pointwise) limit $\lim_n \psi_n(t)$ for every $t \in \mathbb{R}$. It is worth, however, since it points out what we should look for.

From (2.13), (2.50), and (2.53) we get that for any fixed $q \in [0, 1)$

$$\mathbb{E}[X_n] = -\frac{1}{2^{n/2}} + O\left(q^{2^n}\right) \qquad \text{and} \qquad \mathbb{V}[X_n] = \beta(q) + O\left(q^{2^n}\right), \qquad n \to \infty, \tag{2.57}$$

thus

$$\mathbb{E}\left[X_n^2\right] = \beta(q) + \frac{1}{2^n} + O\left(q^{2^n}\right), \qquad n \to \infty. \tag{2.58}$$

One important consequence of (2.58) is [6] that the sequence $F_n$, $n = 0, 1, \ldots$, where $F_n$ is the distribution function of $X_n$, is tight.

**Theorem 3.** For $m = 1, 2, \ldots$ and $q \in [0, 1)$ we have

$$\lim_n \mathbb{E}\left[X_n^m\right] = \beta(q)^{m/2}\mathbb{E}\left[Z^m\right], \tag{2.59}$$

where $\beta(q)$ is given by (2.51), $X_n$ is given by (2.53), and $Z$ is a standard normal random variable. In other words

$$\lim_n \mathbb{E}\left[X_n^{2\ell-1}\right] = 0 \quad \text{and} \quad \lim_n \mathbb{E}\left[X_n^{2\ell}\right] = \frac{(2\ell)!}{2^\ell\,\ell!}\beta(q)^\ell, \qquad \ell \geq 1. \tag{2.60}$$

*Proof.* By differentiating both sides of (2.56) with respect to $t$ we get

$$\begin{aligned}\psi_n'(t) &= \sqrt{2}\, e^{2^{-n/2}\, it}\, \psi_{n-1}\left(\frac{t}{\sqrt{2}}\right)\psi_{n-1}'\left(\frac{t}{\sqrt{2}}\right) + 2^{-n/2} i e^{2^{-n/2}\, it}\, \psi_{n-1}\left(\frac{t}{\sqrt{2}}\right)^2 \\ &\quad + q^{2^{n-1}}\frac{d}{dt}\left[\left(e^{2^{-n/2}\, it} - e^{2\cdot 2^{-n/2}\, it}\right)e^{-2^{n/2}\,\alpha_1(q)it/2}\psi_{n-1}\left(\frac{t}{\sqrt{2}}\right)\right] \\ &\quad + q^{2^n}\frac{d}{dt}\left[\left(e^{2^{-n/2}\, it} - e^{2\cdot 2^{-n/2}\, it}\right)e^{-2^{n/2}\,\alpha_1(q)it}\right].\end{aligned} \tag{2.61}$$

We continue by differentiating $k-1$ times with respect to $t$ both sides of (2.61), where $k \geq 2$. This yields

$$\psi_n^{(k)}(0) = Q_{1n}(0) + Q_{2n}(0) + Q_{3n}(0) + Q_{4n}(0), \tag{2.62}$$

where

$$Q_{1n}(t) := \sqrt{2}\,\frac{d^{k-1}}{dt^{k-1}}\left[e^{2^{-n/2}\,it}\,\psi_{n-1}\left(\frac{t}{\sqrt{2}}\right)\psi'_{n-1}\left(\frac{t}{\sqrt{2}}\right)\right], \tag{2.63}$$

$$Q_{2n}(t) := 2^{-n/2} i\frac{d^{k-1}}{dt^{k-1}}\left[e^{2^{-n/2}\,it}\,\psi_{n-1}\left(\frac{t}{\sqrt{2}}\right)^2\right], \tag{2.64}$$

$$Q_{3n}(t) := q^{2^{n-1}}\frac{d^k}{dt^k}\left[\left(e^{2^{-n/2}\,it} - e^{2\cdot 2^{-n/2}\,it}\right)e^{-2^{n/2}\,\alpha_1(q)it/2}\psi_{n-1}\left(\frac{t}{\sqrt{2}}\right)\right], \tag{2.65}$$

and

$$Q_{4n}(t) := q^{2^n}\frac{d^k}{dt^k}\left[\left(e^{2^{-n/2}\,it} - e^{2\cdot 2^{-n/2}\,it}\right)e^{-2^{n/2}\,\alpha_1(q)it}\right]. \tag{2.66}$$

We will prove the theorem by induction. For $m = 1$ and $m = 2$ the truth of (2.59) follows immediately from (2.57) and (2.58). The inductive hypothesis is that the limit $\lim_n \mathbb{E}\left[X_n^k\right] = i^{-k}\lim_n \psi_n^{(k)}(0)$ satisfies (2.59) or, equivalently, (2.60) for $k = 1, \dots, m-1$. Then, for $k = m$ (where $m \geq 3$):
(i) Formula (2.63) implies

$$Q_{1n}(0) = \frac{1}{(\sqrt{2})^{m-2}}\sum_{j=0}^{m-1}\binom{m-1}{j}\psi_{n-1}^{(m-j)}(0)\,\psi_{n-1}^{(j)}(0) + O\left(\frac{1}{2^{n/2}}\right), \quad n\to\infty; \tag{2.67}$$

(ii) formula (2.64) implies

$$Q_{2n}(0) = O\left(\frac{1}{2^{n/2}}\right), \qquad n\to\infty; \tag{2.68}$$

(iii) since $[e^{2^{-n/2}\,it} - e^{2\cdot 2^{-n/2}\,it}]_{t=0} = 0$, formula (2.65) implies

$$Q_{3n}(0) = O\left(q^{2^{n-1}}2^{mn/2}\right), \qquad n\to\infty; \tag{2.69}$$

(iv) finally, formula (2.66) implies

$$Q_{4n}(0) = O\left(q^{2^n}2^{mn/2}\right), \qquad n\to\infty. \tag{2.70}$$

Thus, for $k = m \geq 3$, by using (2.67), (2.68), (2.69), and (2.70) in (2.62) we can conclude that for every $n \geq 1$ (recall that $\psi_n(0) = 1$)

$$\psi_n^{(m)}(0) = \frac{\psi_{n-1}^{(m)}(0)}{(\sqrt{2})^{m-2}} + \frac{1}{(\sqrt{2})^{m-2}}\sum_{j=1}^{m-1}\binom{m-1}{j}\psi_{n-1}^{(m-j)}(0)\,\psi_{n-1}^{(j)}(0) + \delta_n, \tag{2.71}$$

where $\psi_0^{(m)}(0) = i^m[1-\alpha_1(q)]^m$ and

$$\delta_n = O\left(\frac{1}{2^{n/2}}\right), \qquad n\to\infty. \tag{2.72}$$

Let us set

$$b_n := \frac{1}{(\sqrt{2})^{m-2}} \sum_{j=1}^{m-1} \binom{m-1}{j} \psi_{n-1}^{(m-j)}(0)\, \psi_{n-1}^{(j)}(0) + \delta_n. \tag{2.73}$$

Then, the inductive hypothesis together with (2.72) imply that the limit $\lim_n b_n$ exists in $\mathbb{C}$. We can, therefore, apply Corollary A2 of the Appendix to (2.71) with $z_n = \psi_n^{(m)}(0)$, $\rho = \sqrt{2})^{-(m-2)}$ and $b_n$ as in (2.73) to conclude that the limit $\lim_n \psi_n^{(m)}(0)$ exists in $\mathbb{C}$. This allows us to take limits in (2.71) and obtain

$$\lambda(m) = \frac{\lambda(m)}{(\sqrt{2})^{m-2}} + \frac{1}{(\sqrt{2})^{m-2}} \sum_{j=1}^{m-1} \binom{m-1}{j} \lambda(m-j)\lambda(j)$$

or

$$\left[(\sqrt{2})^{m-2} - 1\right] \lambda(m) = \sum_{j=1}^{m-1} \binom{m-1}{j} \lambda(m-j)\lambda(j) \tag{2.74}$$

where for typographical convenience we have set

$$\lambda(m) := \lim_n \psi_n^{(m)}(0) = i^m \lim_n \mathbb{E}\left[X_n^m\right], \qquad m \geq 1. \tag{2.75}$$

We have, thus, shown that formula (2.74) holds for all $m \geq 1$. If $m = 1$, the sum in the right-hand side of (2.74) is empty, i.e. zero, and, hence we get what we already knew, namely that $\lambda(1) = 0$. If $m = 2$ then (2.74) becomes vacuous ($0 = 0$), but we know that $\lambda(2) = -\beta(q)$. Thus, formula (2.74) gives recursively the values of $\lambda(m)$ for every $m \geq 3$, and by using the induction hypothesis again, namely that $\lambda(k) = i^k \beta(q)^{k/2} \mathbb{E}\left[Z^k\right]$ for $k = 1, \dots, m-1$, it is straightforward to show that $\lambda(m)$, as given by (2.74), equals $i^m \beta(q)^{m/2} \mathbb{E}\left[Z^m\right]$ for every $m$. ∎

A remarkable consequence of Theorem 3 is the following corollary whose proof uses standard arguments, but we decided to include it for the sake of completeness.

**Corollary 2.** Let $\beta(q)$ be the quantity given in (2.51). Then, for any fixed $q \in [0, 1)$ we have

$$X_n \xrightarrow{d} \sqrt{\beta(q)}\, Z \qquad \text{as } n \to \infty, \tag{2.76}$$

where $X_n$ is given by (2.53), $Z$ is a standard normal random variable, and the symbol $\xrightarrow{d}$ denotes convergence in distribution.

*Proof.* As we have seen, the sequence of the distribution functions $F_n(x) = \mathbb{P}\{X_n \leq x\}$, $n = 0, 1, \dots$, is tight.

Suppose that $X_{n_k} \xrightarrow{d} X$, where $X_{n_k}$, $k = 1, 2, \dots$, is a subsequence of $X_n$. Then [6] there is a sequence of random variables $Y_k$, $k = 1, 2, \dots$ converging

a.s. to a random variable $Y$, such that for each $k \geq 1$ the variables $Y_k$ and $X_{n_k}$ have the same distribution (hence the limits $Y$ and $X$ also have the same distribution, since a.s. convergence implies convergence in distribution). Thus, by Theorem 3 we have

$$\lim_k \mathbb{E}\left[Y_k^m\right] = \lim_k \mathbb{E}\left[X_{n_k}^m\right] = \beta(q)^{m/2}\mathbb{E}\left[Z^m\right] \qquad \text{for all integers } m \geq 1. \tag{2.77}$$

In particular, for $m = 2\ell$ we get from (2.77) that $\sup_k \mathbb{E}\left[Y_k^{2\ell}\right] < \infty$ for all $\ell = 1, 2, \ldots$ and, consequently [4], that the sequence $Y_k^m$, $k = 1, 2, \ldots$, is uniformly integrable for every $m \geq 1$. Therefore $Y_k \to Y$ in $L_r(\Omega, \mathcal{F}, \mathbb{P})$, for all $r \geq 0$ (i.e. $\mathbb{E}[|Y_k - Y|^r) \to 0$ as $k \to \infty$, since by $L_r(\Omega, \mathcal{F}, \mathbb{P})$ we denote the space of all real-valued random variables of our probability space $(\Omega, \mathcal{F}, \mathbb{P})$, mentioned in the introduction, with a finite $\|\cdot\|_r$-norm), and (2.77) yields

$$\mathbb{E}\left[X^m\right] = \mathbb{E}\left[Y^m\right] = \lim_k \mathbb{E}\left[Y_k^m\right] = \beta(q)^{m/2}\mathbb{E}\left[Z^m\right] \qquad \text{for all integers } m \geq 1. \tag{2.78}$$

It is well known (and not hard to check) that the moments of the normal distribution satisfy the Carleman's condition and, consequently, a normal distribution is uniquely determined from its moments [4], [6] (alternatively, since the characteristic function of a normal variable is entire, it is uniquely determine by the moments). Therefore, it follows from (2.78) that $X$ and $\sqrt{\beta(q)}\, Z$ have the same distribution. Hence, every subsequence of $X_n$ which converges in distribution, converges to $\sqrt{\beta(q)}\, Z$ and (2.76) is established. ■

From formula (2.53) we have

$$\frac{X_n}{2^{n/2}} = \frac{W_n}{2^n} - \alpha_1(q), \tag{2.79}$$

hence Theorem 3 has the following immediate corollary:

**Corollary 3.** For any $r > 0$ and $q \in [0, 1]$ we have

$$\frac{W_n}{2^n} \to \alpha_1(q) \text{ in } L_r(\Omega, \mathcal{F}, \mathbb{P}), \qquad n \to \infty. \tag{2.80}$$

Finally, let us observe that for any given $\epsilon > 0$ we have (in view of (2.79) and Chebyshev's inequality)

$$\mathbb{P}\left\{\left|\frac{W_n}{2^n} - \alpha_1(q)\right| \geq \epsilon\right\} = \mathbb{P}\left\{\left|\frac{X_n}{2^{n/2}}\right| \geq \epsilon\right\} \leq \frac{1}{2^n \epsilon^2}\mathbb{E}\left[X_n^2\right], \tag{2.81}$$

hence (2.58) yields

$$\mathbb{P}\left\{\left|\frac{W_n}{2^n} - \alpha_1(q)\right| \geq \epsilon\right\} \leq \frac{1}{2^n \epsilon^2}\left[\beta(q) + o(1)\right], \qquad n \to \infty, \tag{2.82}$$

from which it follows by a standard application of the 1st Borel-Cantelli Lemma that if $W_n$, $n = 1, 2, \ldots$, are considered random variables of the space $(\Omega, \mathcal{F}, \mathbb{P})$, then, for any $q \in [0, 1]$

$$\frac{W_n}{2^n} \to \alpha_1(q) \text{ a.s.}, \qquad n \to \infty. \tag{2.83}$$

# 3 The case of a general $N$

Let us now discuss the case of a general $N$, namely the case where $N$ is not necessarily a power of 2. As we have described in the introduction, the first step of the binary search scheme is to test a pool containing samples from all $N$ containers. If this pool is not contaminated, then none of the contents of the $N$ containers is contaminated and we are done. If the pool is contaminated, we form two subpools, one containing samples of the first $\lfloor N/2 \rfloor$ containers and the other containing samples of the remaining $\lceil N/2 \rceil$ containers (recall that $\lfloor N/2 \rfloor + \lceil N/2 \rceil = N$). We continue by testing the first of those subpools. If it is contaminated we split it again into two subpools of $\lfloor \lfloor N/2 \rfloor /2 \rfloor$ and $\lceil \lfloor N/2 \rfloor /2 \rceil$ samples respectively and keep going. We also apply the same procedure to the second subpool of the $\lceil N/2 \rceil$ samples.
Suppose $T(N) = T(N;q)$ is the number of tests required to find all contaminated samples by following the above procedure (thus $T(2^n) = W_n$, where $W_n$ is the random variable studied in the previous section). Then, as in formula (1.2),

$$1 \le T(N) \le 2N - 1. \tag{3.1}$$

In the extreme cases $q = 1$ and $q = 0$ the quantity $T(N)$ becomes deterministic and we have respectively

$$T(N;1) = 1 \qquad \text{and} \qquad T(N;0) = 2N - 1. \tag{3.2}$$

Evidently, $T(N) \le_{st} T(N+1)$, where $\le_{st}$ denotes the usual stochastic ordering (recall that $X \le_{st} Y$ means that $\mathbb{P}\{X > x\} \le \mathbb{P}\{Y > x\}$ for all $x \in \mathbb{R}$). In other words, $T(N)$ is stochastically increasing in $N$. In particular

$$W_{\lfloor \nu \rfloor} \le_{st} T(N) \le_{st} W_{\lceil \nu \rceil}, \qquad \nu := \log_2 N, \tag{3.3}$$

where $\log_2$ is the logarithm to the base 2. Also, it follows easily from a coupling argument that if $q_1 > q_2$, then $T(N;q_1) \le_{st} T(N;q_2)$.

## 3.1 The expectation, the generating function, and the variance of $T(N)$

**Theorem 4.** Let us set

$$\mu(N) = \mu(N;q) := \mathbb{E}[T(N;q)]. \tag{3.4}$$

Then $\mu(N)$ satisfies the recursion

$$\mu(N) = \mu\left(\lfloor N/2 \rfloor\right) + \mu\left(\lceil N/2 \rceil\right) - q^N - q^{\lfloor N/2 \rfloor} + 1, \qquad N \ge 2 \tag{3.5}$$

(of course, $\mu(1) = 1$).

*Proof.* We adapt the proof of Theorem 1. Assume $N \geq 2$ and let $D_N$ be the event that none of the $N$ samples is contaminated. Then

$$\mu(N) = \mathbb{E}[T(N)] = \mathbb{E}[T(N)|D_N]\,\mathbb{P}(D_N) + \mathbb{E}[T(N)|D_N^c]\,\mathbb{P}(D_N^c),$$

and hence

$$\mu(N) = q^N + u(N)\left(1 - q^N\right), \tag{3.6}$$

where for typographical convenience we have set

$$u(N) := \mathbb{E}[T(N)|D_N^c]. \tag{3.7}$$

In order to find a recursive formula for $u(N)$ let us first consider the event $A_N$ that in the group of the $N$ containers none of the first $\lfloor N/2 \rfloor$ contain contaminated samples. Clearly $D_N \subset A_N$ and

$$\mathbb{P}(A_N|D_N^c) = \frac{\mathbb{P}(A_N) - \mathbb{P}(D_N)}{\mathbb{P}(D_N^c)} = \frac{q^{\lfloor N/2 \rfloor} - q^N}{1 - q^N} = \frac{q^{\lfloor N/2 \rfloor}\left(1 - q^{\lceil N/2 \rceil}\right)}{1 - q^N}. \tag{3.8}$$

Likewise, if $B_N$ is the event that in the group of the $N$ containers none of the last $\lceil N/2 \rceil$ contains contaminated samples, then

$$\mathbb{P}(B_N|D_N^c) = \frac{q^{\lceil N/2 \rceil} - q^N}{1 - q^N} = \frac{q^{\lceil N/2 \rceil}\left(1 - q^{\lfloor N/2 \rfloor}\right)}{1 - q^N}. \tag{3.9}$$

Let us also notice that $A_N B_N = D_N$, hence $\mathbb{P}(A_N B_N | D_N^c) = 0$.
Now, (i) given $A_N$ and $D_N^c$ we have that $T(N) = 1 + \tilde{T}(\lceil N/2 \rceil)$, where $\tilde{T}(\lceil N/2 \rceil)$ has the same distribution as $T(\lceil N/2 \rceil)$; (ii) given $B_N$ and $D_N^c$ we have that $T(N) = 2 + T(\lfloor N/2 \rfloor)$; finally (iii) given $(A_N \cup B_N)^c$ and $D_N^c$ we have that $T(N) \overset{d}{=} 1 + T(\lfloor N/2 \rfloor) + \tilde{T}(\lceil N/2 \rceil)$, where, $\tilde{T}(\lceil N/2 \rceil)$ is independent of $T(\lfloor N/2 \rfloor)$ and has the same distribution as $T(\lceil N/2 \rceil)$. Thus (given $D_N^c$ ) by conditioning on the events $A_N$, $B_N$ and $(A_N \cup B_N)^c$ we get, in view of (3.7), (3.8), and (3.9),

$$\begin{aligned} u(N) = {} & \left[1 + u\left(\lceil N/2 \rceil\right)\right] \frac{q^{\lfloor N/2 \rfloor} - q^N}{1 - q^N} + \left[2 + u\left(\lfloor N/2 \rfloor\right)\right] \frac{q^{\lceil N/2 \rceil} - q^N}{1 - q^N} \\ & + \left[1 + u\left(\lceil N/2 \rceil\right) + u\left(\lfloor N/2 \rfloor\right)\right] \left(1 - \frac{q^{\lfloor N/2 \rfloor} - q^N}{1 - q^N} - \frac{q^{\lceil N/2 \rceil} - q^N}{1 - q^N}\right) \end{aligned}$$

or

$$\begin{aligned} \left(1 - q^N\right) u(N) = {} & \left(1 - q^{\lfloor N/2 \rfloor}\right) u\left(\lfloor N/2 \rfloor\right) + \left(1 - q^{\lceil N/2 \rceil}\right) u\left(\lceil N/2 \rceil\right) \\ & - 2q^N + q^{\lceil N/2 \rceil} + 1 \end{aligned} \tag{3.10}$$

from which, in view of (3.6), formula (3.5) follows immediately. ■

In the case where $N = 2^n$ formula (3.5) reduces to (2.10).

**Remark 2.** By easy induction formula (3.5) implies that, for $N \geq 2$ the quantity $\mu(N;q) = \mathbb{E}[T(N;q)]$ is a polynomial in $q$ of degree $N$ whose coefficients are $\leq 0$, except for the constant term which is equal to $2N-1$. The leading term of $\mu(N;q)$ is $-q^N$. In the special case where $N = 2^n$ these properties of $\mu(N;q)$ follow trivially from (2.3).

For instance,

$$\begin{aligned} &\mu(3) = 5 - 2q - q^2 - q^3, \qquad \mu(5) = 9 - 3q - 3q^2 - q^3 - q^5, \\ &\text{and} \\ &\mu(1000) = 1999 - 512q - 720q^2 - 48q^3 - 336q^4 - 48q^7 - 144q^8 - 48q^{15} \\ &\quad - 48q^{16} - 40q^{31} - 8q^{32} - 16q^{62} - 8q^{63} - 12q^{125} - 6q^{250} - 3q^{500} - q^{1000}, \end{aligned} \tag{3.11}$$

while $\mu(2)$, $\mu(4)$, and $\mu(1024)$ are given by (2.12), since $W_1 = T(2)$, $W_2 = T(4)$, and $W_{10} = T(1024)$.
With the help of (3.5) one can obtain an extension of formula (2.3) valid for any $N$. We first need to introduce some convenient notation:

$$\iota(N) := \lfloor N/2 \rfloor, \qquad \iota^k(N) = (\underbrace{\iota \circ \cdots \circ \iota}_{k \text{ iterates}})(N), \tag{3.12}$$

so that $\iota^1(N) = \iota(N)$, while we also have the standard convention $\iota^0(N) = N$ (we furthermore have $\iota^{-1}(N) = \{2N, 2N+1\}$). For example, if $N \geq 2$ and, as in (3.3), $\nu = \log_2 N$, then

$$\iota^{\lfloor \nu \rfloor}(N) = 1, \qquad \text{while} \qquad \iota^{\lfloor \nu \rfloor - 1}(N) = 2 \text{ or } 3 \quad \text{and} \quad \iota^{\lfloor \nu \rfloor + 1}(N) = 0. \tag{3.13}$$

**Corollary 4.** Let $\mu(N)$ be as in (3.4) and

$$\begin{aligned} \varepsilon_\mu(N) = \varepsilon_\mu(N;q) &:= q^{\lfloor N/2 \rfloor} + q^N - q^{\lfloor (N+1)/2 \rfloor} - q^{N+1} \\ &= q^{\lfloor N/2 \rfloor} - q^{\lceil N/2 \rceil} + q^N - q^{N+1} \end{aligned} \tag{3.14}$$

(the last equality follows from the fact that $\lfloor (N+1)/2 \rfloor = \lceil N/2 \rceil$). Then

$$\mu(N) = 2N - 1 - (N-1)q - (N-1)q^2 + \sum_{n=2}^{N-1} \sum_{k=1}^{\lfloor \log_2 n \rfloor} \varepsilon_\mu\left(\iota^{k-1}(n)\right), \qquad N \geq 1 \tag{3.15}$$

(if $N = 1$ or $N = 2$, then the double sum in the right-hand side is 0).

*Proof.* By setting

$$\Delta_\mu(N) = \Delta_\mu(N;q) := \mu(N+1;q) - \mu(N;q) \tag{3.16}$$

and by recalling that $\lfloor (N+1)/2 \rfloor = \lceil N/2 \rceil$ and $\lceil (N+1)/2 \rceil = \lfloor N/2 \rfloor + 1$, formula (3.5) implies (in view of (3.12) and (3.14))

$$\Delta_\mu(N) = \Delta_\mu\left(\lfloor N/2 \rfloor\right) + \varepsilon_\mu(N) = \Delta_\mu\big(\iota(N)\big) + \varepsilon_\mu(N), \qquad N \geq 2. \quad (3.17)$$

From (3.17) we obtain

$$\sum_{k=1}^{\lfloor \log_2 N \rfloor} \left[\Delta_\mu\left(\iota^{k-1}(N)\right) - \Delta_\mu\left(\iota^{k}(N)\right)\right] = \sum_{k=1}^{\lfloor \log_2 N \rfloor} \varepsilon_\mu\left(\iota^{k-1}(N)\right), \qquad N \geq 2,$$

or, in view of (3.13),

$$\Delta_\mu(N) - \Delta_\mu(1) = \sum_{k=1}^{\lfloor \log_2 N \rfloor} \varepsilon_\mu\left(\iota^{k-1}(N)\right), \qquad N \geq 1, \qquad (3.18)$$

where $\Delta_\mu(1) = \mu(2) - \mu(1) = \mathbb{E}[T(2)] - \mathbb{E}[T(1)] = \mathbb{E}[W_1] - \mathbb{E}[W_0] = 2 - q - q^2$ (in the case $N = 1$, formula (3.18) is trivially true, since the sum in the right-hand side is empty). Consequently, (3.18) becomes

$$\Delta_\mu(N) = \mu(N+1) - \mu(N) = 2 - q - q^2 + \sum_{k=1}^{\lfloor \log_2 N \rfloor} \varepsilon_\mu\left(\iota^{k-1}(N)\right), \qquad N \geq 1. \qquad (3.19)$$

from which (3.15) follows immediately. ■

With the help of (3.15) we have created the following table for the average-case aspect ratio:

**Table 2:** *Values of $\mathbb{E}[T(N)]/N$ for certain choices of $N \neq 2^n$ and $p$.*

| | $p = .005$ | $p = .01$ | $p = .05$ | $p = .1$ | $p = .15$ | $p = .2$ |
|---|---|---|---|---|---|---|
| $N = 3$ | .345 | .357 | .447 | .554 | .654 | .749 |
| $N = 5$ | .217 | .234 | .362 | .510 | .645 | .768 |
| $N = 6$ | .186 | .205 | .348 | .510 | .656 | .787 |
| $N = 7$ | .163 | .184 | .337 | .510 | .663 | .800 |
| $N = 12$ | .110 | .136 | .325 | .526 | .696 | .843 |
| $N = 48$ | .061 | .099 | .345 | .571 | .751 | .902 |
| $N = 96$ | .057 | .099 | .354 | .582 | .762 | .912 |
| $N = 200$ | .057 | .101 | .358 | .585 | .765 | .917 |
| $N = 389$ | .058 | .104 | .362 | .589 | .769 | .920 |
| $N = 768$ | .059 | .105 | .363 | .591 | .771 | .922 |
| $N = 1000$ | .060 | .106 | .365 | .594 | .776 | .929 |

(as in Table 1, if we compare the above values with the graphs, given in [1] and [2], of the optimal average-case aspect ratio as a function of $p$, we can see that, if $p \leq 0.15$, our procedure is near-optimal).

**Remark 3.** Let us look at the case $N = 2^n M$, where $M$ is a fixed odd integer. We set

$$\tau_n := \mathbb{E}\left[T(N)\right] = \mathbb{E}\left[T(2^n M)\right], \tag{3.20}$$

so that

$$\tau_0 = \mathbb{E}\left[T(M)\right] = \mu(M). \tag{3.21}$$

Then (3.5) becomes

$$\tau_n = 2\tau_{n-1} - q^{2^n M} - q^{2^{n-1} M} + 1, \qquad n \geq 1. \tag{3.22}$$

From (3.21) and (3.22) it follows that

$$\begin{aligned} \frac{\tau_n}{2^n} &= \mu(M) + \sum_{k=1}^{n} \frac{1 - q^{2^k M} - q^{2^{k-1} M}}{2^k} \\ &= \mu(M) + 1 - \frac{1}{2^n} - \sum_{k=1}^{n} \frac{q^{2^k M} + q^{2^{k-1} M}}{2^k} \qquad n \geq 0, \end{aligned} \tag{3.23}$$

thus

$$\frac{\mathbb{E}\left[T(2^n M)\right]}{2^n M} = \frac{\tau_n}{2^n M} = \frac{\mu(M) + 1}{M} - \frac{1}{2^n M} - \sum_{k=1}^{n} \frac{q^{2^k M} + q^{2^{k-1} M}}{2^k M} \tag{3.24}$$

for all $n \geq 0$. It follows that, as $n \to \infty$,

$$\frac{\mathbb{E}\left[T(2^n M)\right]}{2^n M} \to \alpha_M(q) \qquad \text{uniformly in } \ q \in [0, 1], \tag{3.25}$$

where

$$\alpha_M(q) := \frac{\mu(M)+1}{M} - \sum_{k=1}^{\infty} \frac{q^{2^k M} + q^{2^{k-1}M}}{2^k M}$$
$$= \frac{2\mu(M)+2-q^M}{2M} - \frac{3}{2}\sum_{k=1}^{\infty} \frac{q^{2^k M}}{2^k M} \tag{3.26}$$

(as in the case of Remark 1, the convergence in (3.25) is much stronger, since for every $m = 1, 2, \dots$ the $m$-th derivative of $\mathbb{E}[T(2^n M)]/(2^n M)$ with respect to $q$ converges to the $m$-th derivative of $\alpha_M(q)$ uniformly on $[0, 1-\epsilon]$ for any $\epsilon > 0$).
From (3.26) it is obvious that if $M_1 \neq M_2$, then $\alpha_{M_1}(q) \neq \alpha_{M_2}(q)$, except for at most countably many values of $q$ (notice, e.g., that $\alpha_M(0) = 2$ and $\alpha_M(1) = 0$ for all $M$). Therefore, it follows from (3.25) that for $q \in (0, 1)$

$$\text{the limit} \qquad \lim_N \frac{\mathbb{E}\left[T(N;q)\right]}{N} \qquad \text{does not exist.} \tag{3.27}$$

**Comment.** The reason we consider the case $N = 2^n M$, where $M$ is odd, is mainly theoretical: We use these values of $N$ and the corresponding subsequences in order to demonstrate the non-convergence of $T(N)/N$. We do not claim any practical value of this choice of $N$.

**Open Question.** Is it true that $\alpha_1(q) = \limsup_N \mathbb{E}[T(N;q)]/N$?

Next, we consider the generating function of $T(N)$. Using the approach of the proofs of Theorems 2 and 4 we can derive the following result:

**Theorem 5.** Let

$$g(z;N) = g(z;N;q) := \mathbb{E}\left[z^{T(N;q)}\right]. \tag{3.28}$$

Then $g(z;N)$ satisfies the recursion

$$g(z;N) = zg\left(z;\lfloor N/2\rfloor\right) g\left(z;\lceil N/2\rceil\right) + (z-z^2)q^{\lfloor N/2\rfloor} g\left(z;\lceil N/2\rceil\right) + (z-z^2)q^N \tag{3.29}$$

for $N \geq 2$ (of course, $g(z;1) = z$).

Notice that in the case where $N = 2^n$, formula (3.29) reduces to (2.31).
By setting $z = e^{it}$ in (3.29) it follows that the characteristic function

$$\phi(t;N) = \phi(t;N;q) := \mathbb{E}\left[e^{itT(N;q)}\right] \tag{3.30}$$

of $T(N) = T(N;q)$ satisfies the recursion

$$\phi(t;N) = e^{it}\phi\left(t;\lfloor N/2\rfloor\right)\phi\left(t;\lceil N/2\rceil\right) + \left(e^{it} - e^{2it}\right) q^{\lfloor N/2\rfloor}\phi\left(t;\lceil N/2\rceil\right)$$
$$+ \left(e^{it} - e^{2it}\right) q^N, \qquad N \geq 2, \tag{3.31}$$

with, of course, $\phi(t;1) = e^{it}$. For example, for $N = 2$ formula (3.31) confirms the value of $\phi(t;2)$ given in (2.38).
By differentiating (3.29) twice with respect to $z$ and then setting $z = 1$ and using (3.5) we get the following corollary:

**Corollary 5.** The variance

$$\sigma^2(N) = \sigma^2(N;q) := \mathbb{V}\left[T(N;q)\right] \tag{3.32}$$

satisfies the recursion

$$\begin{aligned}\sigma^2(N) = \sigma^2\left(\lfloor N/2 \rfloor\right) + \sigma^2\left(\lceil N/2 \rceil\right) + 2\mu(N)q^N + 2\mu\left(\lfloor N/2 \rfloor\right)q^{\lfloor N/2 \rfloor} \\ + q^{2N} - q^{2\lfloor N/2 \rfloor} - 3q^N - q^{\lfloor N/2 \rfloor}, \qquad N \geq 2\end{aligned} \tag{3.33}$$

(of course, $\sigma^2(1) = 0$).

In the case $N = 2^n$ formula (3.33) reduces to (2.44) (e.g., if $N = 2$, then (3.33) agrees with (2.46)).
As we have mentioned, in the extreme cases where $q = 0$ or $q = 1$ the variable $T(N)$ becomes deterministic and, consequently, $\sigma^2(N;0) = \sigma^2(N;1) = 0$, which is in agreement with (3.33).

**Corollary 6.** For $q \in (0,1)$ there are constants $0 < c_1 < c_2$, depending on $q$, such that $\sigma^2(N) = \mathbb{V}\left[T(N)\right]$ satisfies

$$c_1 N \leq \sigma^2(N) \leq c_2 N \qquad \text{for all } \; N \geq 2. \tag{3.34}$$

*Proof.* For $N \geq 2$ we have $\mathbb{P}\{T(N) = 1\} = q^N$, which implies that $\mu(N) \geq q^N + 2(1 - q^N) = 2 - q^N$. Hence

$$\begin{aligned}&2\mu(N)q^N + 2\mu\left(\lfloor N/2 \rfloor\right)q^{\lfloor N/2 \rfloor} + q^{2N} - q^{2\lfloor N/2 \rfloor} - 3q^N - q^{\lfloor N/2 \rfloor} \\ &\geq 2\left(2 - q^N\right)q^N + 2\left(2 - q^{\lfloor N/2 \rfloor}\right)q^{\lfloor N/2 \rfloor} + q^{2N} - q^{2\lfloor N/2 \rfloor} - 3q^N - q^{\lfloor N/2 \rfloor} \\ &= q^N - q^{2N} + 3q^{\lfloor N/2 \rfloor} - 3q^{2\lfloor N/2 \rfloor} > 0.\end{aligned} \tag{3.35}$$

Using (3.35) in (3.33) implies

$$\sigma^2(N) > \sigma^2\left(\lfloor N/2 \rfloor\right) + \sigma^2\left(\lceil N/2 \rceil\right), \qquad N \geq 2. \tag{3.36}$$

Now, as we have already seen, $\sigma^2(1;q) \equiv 0$, while by (2.46) we have $\sigma^2(2;q) = q(1-q)(q^2 + 3q + 1)$. We can, thus, choose

$$c_1 = \frac{\sigma^2(2;q)}{3} = \frac{q(1-q)(q^2+3q+1)}{3} > 0. \tag{3.37}$$

Then, (3.37) implies that $\sigma^2(2) = 3c_1 > 2c_2$ and from (3.36) we have $\sigma^2(3) > \sigma^2(2) = 3c_1$, i.e. the first inequality in (3.34) is valid for $N = 2$ and $N = 3$.

Finally, the inequality $c_1 N \leq \sigma^2(N)$ for every $N \geq 2$ follows easily from (3.36) by induction.
To establish the second inequality of (3.34) we set

$$\Delta_\sigma(N) := \sigma^2(N+1) - \sigma^2(N). \tag{3.38}$$

Then, formula (3.33) implies (in view of (3.12))

$$\Delta_\sigma(N) = \Delta_\sigma\left(\lfloor N/2 \rfloor\right) + \varepsilon_\sigma(N) = \Delta_\sigma\big(\iota(N)\big) + \varepsilon_\sigma(N), \qquad N \geq 2, \tag{3.39}$$

where

$$\begin{aligned}\varepsilon_\sigma(N) = \varepsilon_\sigma(N;q) := {} & 2\left[\mu(N+1)q - \mu(N)\right] q^N + 2\mu\left(\lceil N/2 \rceil\right) q^{\lceil N/2 \rceil} \\ & - 2\mu\left(\lfloor N/2 \rfloor\right) q^{\lfloor N/2 \rfloor} + q^{\lfloor N/2 \rfloor} - q^{\lceil N/2 \rceil} + q^{2\lfloor N/2 \rfloor} \\ & - q^{2\lceil N/2 \rceil} + \left(3 - q^N - q^{N+1}\right) q^N (1-q). \end{aligned} \tag{3.40}$$

Observe that by (3.40) and (3.1) we have that

$$\varepsilon_\sigma(N) = O\left(N q^{\lfloor N/2 \rfloor}\right), \qquad N \to \infty. \tag{3.41}$$

Now, as in the case of (3.18) we have, in view of (3.39),

$$\Delta_\sigma(N) = \Delta_\sigma(1) + \sum_{k=1}^{\lfloor \log_2 N \rfloor} \varepsilon_\sigma\left(\iota^{k-1}(N)\right), \qquad N \geq 2, \tag{3.42}$$

where $\Delta_\sigma(1) = \sigma^2(2) - \sigma^2(1) = q(1-q)(q^2+3q+1)$. From (3.42) and (3.41) we get that $\Delta_\sigma(N)$ is bounded. Therefore, the second inequality of (3.34) follows immediately from (3.38). ■

Of course, in the case where $N = 2^n$ we have formula (2.50), which is much more precise than (3.34).
Let us also notice that with the help of (3.38), (3.42), and (3.15) we can get a messy, yet explicit formula for $\sigma^2(N)$, extending (2.40) to the case of a general $N$.

## 3.2 The behavior of $T(N)$ as $N \to \infty$

We start with a Lemma.

**Lemma 1.** Let $\sigma^2(N) = \sigma^2(N;q) = \mathbb{V}\left[T(N;q)\right]$ as in the previous subsection. Then for any fixed $q \in (0,1)$ we have

$$\frac{\sigma^2\left(\lfloor N/2 \rfloor\right)}{\sigma^2(N)} = \frac{1}{2} + O\left(\frac{1}{N}\right) \quad \text{and} \quad \frac{\sigma^2\left(\lceil N/2 \rceil\right)}{\sigma^2(N)} = \frac{1}{2} + O\left(\frac{1}{N}\right), \quad N \to \infty. \tag{3.43}$$

*Proof.* As we have seen in the proof of Corollary 6, $\Delta_\sigma(N) = \sigma^2(N+1) - \sigma^2(N)$ is bounded, hence

$$\sigma^2\left(\lceil N/2 \rceil\right) = \sigma^2\left(\lfloor N/2 \rfloor\right) + O(1), \qquad N \to \infty. \tag{3.44}$$

Now, if we divide (3.33) by $\sigma^2(N)$ and invoke (3.34) we get

$$1 = \frac{\sigma^2\left(\lfloor N/2 \rfloor\right)}{\sigma^2(N)} + \frac{\sigma^2\left(\lceil N/2 \rceil\right)}{\sigma^2(N)} + O\left(\frac{q^{N/2}}{N}\right), \qquad N \to \infty. \tag{3.45}$$

Therefore, by using (3.44) in (3.45) and invoking again (3.34) we obtain

$$1 = \frac{2\sigma^2\left(\lfloor N/2 \rfloor\right)}{\sigma^2(N)} + O\left(\frac{1}{N}\right), \qquad N \to \infty, \tag{3.46}$$

which is equivalent to the first equality of (3.43). The second equality follows immediately. ∎

In order to determine the limiting behavior of $T(N)$, as $N \to \infty$, for any fixed $q \in (0,1)$, it is natural to work with the normalized variables

$$Y(N) := \frac{T(N) - \mu(N)}{\sigma(N)}, \qquad N = 2, 3, \dots \tag{3.47}$$

($Y(1)$ is not defined). Obviously,

$$\mathbb{E}[Y(N)] = 0 \quad \text{and} \quad \mathbb{V}[Y(N)] = \mathbb{E}[Y(N)^2] = 1 \qquad \text{for all } N \geq 2, \tag{3.48}$$

thus [6] the sequence of the distribution functions of $Y(N)$ is tight.
Let us, also, introduce the characteristic functions

$$\psi(t; N) := \mathbb{E}\left[e^{itY(N)}\right], \qquad N = 2, 3, \dots . \tag{3.49}$$

In view of (3.47) and (3.30) we have

$$\psi(t; N) = \phi\left(\frac{t}{\sigma(N)}; N\right) e^{-it\mu(N)/\sigma(N)t} \tag{3.50}$$

and by using (3.50) in (3.31) and then invoking (3.5) we obtain the recursion

$$\begin{aligned}\psi(t; N) = {} & e^{\frac{q^N + q^{\lfloor N/2 \rfloor}}{\sigma(N)} it} \psi\left(\frac{\sigma(\lfloor N/2 \rfloor)}{\sigma(N)} t; \lfloor N/2 \rfloor\right) \psi\left(\frac{\sigma(\lceil N/2 \rceil)}{\sigma(N)} t; \lceil N/2 \rceil\right) \\ & + q^{\lfloor N/2 \rfloor}\left(e^{\frac{1}{\sigma(N)} it} - e^{\frac{2}{\sigma(N)} it}\right) e^{\frac{\mu(\lceil N/2 \rceil) - \mu(N)}{\sigma(N)} it} \psi\left(\frac{\sigma(\lceil N/2 \rceil)}{\sigma(N)} t; \lceil N/2 \rceil\right) \\ & + q^N\left(e^{\frac{1}{\sigma(N)} it} - e^{\frac{2}{\sigma(N)} it}\right) e^{-\frac{\mu(N)}{\sigma(N)} it}, \qquad N \geq 4, \end{aligned} \tag{3.51}$$

with $\psi(t;2)$ and $\psi(t;3)$ taken from (3.49). Actually, if we define $\psi(0;1):=1$, then (3.51) is valid for $N\geq 2$.
We are now ready to present a general result, which can be viewed as an extension of Theorem 3 to the case of an arbitrary $N$.

**Theorem 6.** For $m=1,2,\dots$ and $q\in(0,1)$ we have

$$\lim_N \mathbb{E}\left[Y(N)^m\right]=\mathbb{E}\left[Z^m\right], \tag{3.52}$$

where $Y(N)$ is given by (3.47) and $Z$ is a standard normal random variable. In other words

$$\lim_N \mathbb{E}\left[Y(N)^{2\ell-1}\right]=0 \quad \text{and} \quad \lim_N \mathbb{E}\left[Y(N)^{2\ell}\right]=\frac{(2\ell)!}{2^\ell\,\ell!}, \qquad \ell\geq 1. \tag{3.53}$$

*Proof.* We will follow the approach of the proof of Theorem 3.
We apply induction. By (3.48) we know that (3.53) is valid in the special cases $m=1$ and $m=2$. The inductive hypothesis is that the limit $\lim_N \mathbb{E}\left[Y(N)^k\right]=i^{-k}\lim_N \psi^{(k)}(0;N)$ (where $\psi^{(k)}(t;N)$ denotes the $k$-th derivative of $\psi(t;N)$ with respect to $t$) satisfies (3.53) for $k=1,\dots,m-1$. Then, for $k=m$ (where $m\geq 3$) formula (3.51) implies

$$\begin{aligned}\psi^{(m)}(0;N)&=\frac{\sigma^m\left(\lfloor N/2\rfloor\right)}{\sigma^m(N)}\,\psi^{(m)}(0;\lfloor N/2\rfloor)+\frac{\sigma^m\left(\lceil N/2\rceil\right)}{\sigma^m(N)}\,\psi^{(m)}(0;\lceil N/2\rceil)\\&\quad+\sum_{j=1}^{m-1}\binom{m}{j}\frac{\sigma^j\left(\lfloor N/2\rfloor\right)\sigma^{m-j}\left(\lceil N/2\rceil\right)}{\sigma^m(N)}\,\psi^{(j)}(0;\lfloor N/2\rfloor)\psi^{(m-j)}(0;\lceil N/2\rceil)\\&\quad+\delta(N),\end{aligned} \tag{3.54}$$

where

$$\delta(N)=O\left(q^{N/2}\right), \qquad N\to\infty. \tag{3.55}$$

Let us set

$$\begin{aligned}b(N):=&\sum_{j=1}^{m-1}\binom{m}{j}\frac{\sigma^j\left(\lfloor N/2\rfloor\right)\sigma^{m-j}\left(\lceil N/2\rceil\right)}{\sigma^m(N)}\,\psi^{(j)}(0;\lfloor N/2\rfloor)\psi^{(m-j)}(0;\lceil N/2\rceil)\\&+\delta(N).\end{aligned} \tag{3.56}$$

Then, the inductive hypothesis together with Lemma 1 and (3.55) imply that the limit $\lim_N b(N)$ exists in $\mathbb{C}$. We can, therefore, apply Corollary A1 of the Appendix to (3.54) with $z(N)=\psi^{(m)}(0;N)$, $\rho_1(N)=\sigma^m\left(\lfloor N/2\rfloor\right)/\sigma^m(N)$, $\rho_2(N)=\sigma^m\left(\lceil N/2\rceil\right)/\sigma^m(N)$ (so that, in view of Lemma 1 and (4.15), $\rho_1=\rho_2=1/2^{m/2}$), and $b(N)$ as in (3.56), to conclude that

$$\text{the limit} \quad \Lambda(m):=\lim_N \psi^{(m)}(0;N) \quad \text{exists in } \mathbb{C}. \tag{3.57}$$

Thus by taking limits in (3.54) we obtain (in view of (3.57) and Lemma 1)

$$\Lambda(m) = \frac{\Lambda(m)}{2^{m/2}} + \frac{\Lambda(m)}{2^{m/2}} + \frac{1}{2^{m/2}} \sum_{j=1}^{m-1} \binom{m}{j} \Lambda(m-j)\Lambda(j)$$

or

$$\left(2^{m/2} - 2\right) \Lambda(m) = \sum_{j=1}^{m-1} \binom{m}{j} \Lambda(m-j)\Lambda(j). \tag{3.58}$$

We have, thus, shown that formula (3.58) holds for all $m \geq 1$. If $m = 1$, the sum in the right-hand side of (3.58) is empty, i.e. zero, and, hence we get what we already knew, namely that $\Lambda(1) = 0$. If $m = 2$ then (3.58) becomes vacuous ($0 = 0$), but we know that $\Lambda(2) = -1$. Thus, formula (3.58) gives recursively the values of $\Lambda(m)$ for every $m \geq 3$, and by using the induction hypothesis again, namely that $\Lambda(k) = i^k \mathbb{E}\left[Z^k\right]$ for $k = 1, \dots, m-1$, it is straightforward to show that $\Lambda(m)$, as given by (3.58), equals $i^m \mathbb{E}\left[Z^m\right]$ for every $m$. ■

Theorem 6 together with the fact that the sequence $Y(N)$ is tight imply the following corollary which can be considered an extension of Corollary 2. Its proof is omitted since it is just a repetition of the proof of Corollary 2.

**Corollary 7.** For any fixed $q \in (0, 1)$ we have

$$Y(N) = \frac{T(N) - \mathbb{E}[T(N)]}{\sqrt{\mathbb{V}[T(N)]}} \xrightarrow{d} Z \qquad \text{as } N \to \infty, \tag{3.59}$$

where $Y(N)$ is given by (3.47) and $Z$ is a standard normal random variable.

Let us, also, mention that, since $\sigma(N) = \sqrt{\mathbb{V}[T(N)]} = O(\sqrt{N})$, a rather trivial consequence of Theorem 6 is that for any given $\epsilon > 0$ we have that

$$\frac{T(N) - \mathbb{E}[T(N)]}{N^{\frac{1}{2}+\epsilon}} \to 0 \text{ in } L_r(\Omega, \mathcal{F}, \mathbb{P}), \qquad N \to \infty, \tag{3.60}$$

for any $r > 0$ and any $q \in [0, 1]$.

It follows that if $N_n$, $n = 1, 2, \dots$, is a sequence of integers (with $\lim_n N_n = \infty$) such that the limit $\lim_n \mathbb{E}[T(N_n)]/N_n$ exists (e.g., $N_n = 2^n$ or $N_n = 3 \cdot 2^n$), then

$$\frac{T(N_n)}{N_n} \to \lim_n \frac{\mathbb{E}[T(N_n)]}{N_n}, \qquad n \to \infty, \tag{3.61}$$

in the $L_r(\Omega, \mathcal{F}, \mathbb{P})$-sense, for every $r > 0$ (the case $N_n = 2^n$ is treated by Corollary 3).

Also, by Chebyshev's inequality (applied to a sufficiently high power of $\left(T(N) - \mathbb{E}[T(N)]\right)/N^{\frac{1}{2}+\epsilon}$) and the 1st Borel-Cantelli Lemma it follows that

$$\frac{T(N) - \mathbb{E}[T(N)]}{N^{\frac{1}{2}+\epsilon}} \to 0 \text{ a.s.}, \qquad N \to \infty, \tag{3.62}$$

for any $q \in [0,1]$ and any given $\epsilon > 0$. In this case, again, if $\lim_n \mathbb{E}[T(N_n)]/N_n$ exists for some sequence $N_n \to \infty$, then

$$\frac{T(N_n)}{N_n} \to \lim_n \frac{\mathbb{E}[T(N_n)]}{N_n} \quad \text{a.s.}, \qquad n \to \infty. \tag{3.63}$$

Finally, let us mention a query: As we have already said in the end of Subsection 2.1, for a given $q = 1 - p$ there is an optimal value $m = m(q)$ of $N$ which minimizes $\mathbb{E}[T(N)]/N$. Tables 1 and 2 suggest that $m$ increases as $p$ decreases and, as we have seen in the last paragraph of Subsection 2.1, a first guess is that $m \sim c/p$ as $p \to 0^+$, for some constant $c > 0$. It will be, then, interesting to understand the limiting behavior of $T(\lfloor c/p \rfloor)$ as $p \to 0^+$, for a given $c > 0$. Actually, in this case we do not believe that the limiting distribution will be normal.

# 4 Appendix

This section contains some auxiliary material. First we present a proposition which is used in the justification of formula (2.17).

**Proposition A1.** We have

$$\int_1^\infty q^{2^x} dx = -\frac{\ln(1-q)}{\ln 2} + O(1), \qquad q \to 1^-. \tag{4.1}$$

*Proof.* By applying the substitution $\xi = |\ln q| 2^x$ we get

$$\int_1^\infty q^{2^x} dx = \int_1^\infty e^{(\ln q)2^x} dx = \frac{1}{\ln 2} \int_{2|\ln q|}^\infty \frac{e^{-\xi}}{\xi} d\xi. \tag{4.2}$$

Since we are interested in the limiting behavior of the integral as $q \to 1^-$, we can assume that $1/\sqrt{e} < q < 1$ (i.e. $1 > 2\,|\ln q| > 0$). Then, (4.2) implies

$$\begin{aligned} \int_1^\infty q^{2^x} dx &= \frac{1}{\ln 2} \int_{2|\ln q|}^1 \frac{d\xi}{\xi} + \frac{1}{\ln 2} \int_{2|\ln q|}^1 \frac{e^{-\xi} - 1}{\xi} d\xi + \frac{1}{\ln 2} \int_1^\infty \frac{e^{-\xi}}{\xi} d\xi \\ &= -\frac{\ln|\ln q|}{\ln 2} + O(1), \qquad q \to 1^- \end{aligned}$$

and (4.1) follows from the fact that $\ln|\ln q| = \ln(1-q) + o(1)$ as $q \to 1^-$. ■

Next, we give the proof equation (2.52).

*Proof of* (2.52). Equation (2.51) can be written as

$$\beta(q) = \sum_{k=1}^\infty f_k(q) + \sum_{k=1}^\infty \frac{q^{2^{k+1}} + q^{3 \cdot 2^{k-1}} - 5q^{2^k} - 3q^{2^{k-1}}}{2^k}, \tag{4.3}$$

where

$$f_k(q) := \left(2 - \sum_{j=1}^{k} \frac{q^{2^j} + q^{2^{j-1}}}{2^j}\right)\left(2q^{2^k} + q^{2^{k-1}}\right) \tag{4.4}$$

(in particular, $f_k(1) = 3(2 - \sum_{j=1}^{k} 2 \cdot 2^{-j}) = 6 \cdot 2^{-k}$). Taking limits in (4.3) we obtain

$$\begin{aligned}\lim_{q \to 1^-} \beta(q) &= \lim_{q \to 1^-} \sum_{k=1}^{\infty} f_k(q) + \lim_{q \to 1^-} \sum_{k=1}^{\infty} \frac{q^{2^{k+1}} + q^{3 \cdot 2^{k-1}} - 5q^{2^k} - 3q^{2^{k-1}}}{2^k} \\ &= \lim_{q \to 1^-} \sum_{k=1}^{\infty} f_k(q) - 6 \sum_{k=1}^{\infty} \frac{1}{2^k} = \lim_{q \to 1^-} \sum_{k=1}^{\infty} f_k(q) - 6. \end{aligned} \tag{4.5}$$

To finish the proof of (2.52) we need to justify the passing of the limit inside the summation in the last term of (4.5). We can do that by invoking the Dominated Convergence Theorem. Thus, we need to dominate $f_k(q)$ by a summable bound. It is obvious from (4.4) that $f_k(q) \geq 0$, hence we only need to find a good upper bound for $f_k(q)$.
In view of (4.4) we have (since $0 \leq q \leq 1$)

$$f_k(q) \leq 6F_k(q), \qquad \text{where} \qquad F_k(q) := \left(1 - \sum_{j=1}^{k} \frac{q^{2^j}}{2^j}\right) q^{2^{k-1}}. \tag{4.6}$$

The derivative of $F_k(q)$ is

$$\begin{aligned} F_k'(q) &= 2^{k-1} q^{2^{k-1}-1} \left(1 - \sum_{j=1}^{k} \frac{q^{2^j}}{2^j}\right) - q^{2^{k-1}} \sum_{j=1}^{k} q^{2^j - 1} \\ &= 2^{k-1} q^{2^{k-1}-1} \left(1 - \sum_{j=1}^{k} \frac{q^{2^j}}{2^j} - \frac{1}{2^{k-1}} \sum_{j=1}^{k} q^{2^j}\right). \end{aligned} \tag{4.7}$$

It is clear from (4.7) that $F_k'(q)$ has a unique zero, say $q_k^*$, in $(0,1)$ and $F_k'(q) > 0$ for $0 < q < q_k^*$, while $F_k'(q) < 0$ for $q_k^* < q < 1$. Also, in view of (4.7), $F_k'(q_k^*) = 0$ is equivalent to the equality

$$1 - \sum_{j=1}^{k} \frac{(q_k^*)^{2^j}}{2^j} = \frac{1}{2^{k-1}} \sum_{j=1}^{k} (q_k^*)^{2^j}. \tag{4.8}$$

Now, $F_k(q)$ attains its maximum on $[0,1]$ at $q = q_k^*$. Hence, in view of (4.6) and (4.8) we get

$$F_k(q) \leq F_k(q_k^*) = \frac{(q_k^*)^{2^{k-1}}}{2^{k-1}} \sum_{j=1}^{k} (q_k^*)^{2^j} \leq \frac{k}{2^{k-1}}, \qquad q \in [0,1]. \tag{4.9}$$

Therefore, as it follows from (4.6) and (4.9), $f_k(q)$ is bounded by $3k2^{-k}$ (which is, obviously, summable) and, consequently, dominated convergence applies. ■

Finally, we present a lemma with two corollaries which are used in the proofs of some results of Sections 2 and 3.

**Lemma A1.** Let $\varepsilon(N)$, $N = 1, 2, \ldots$, be a sequence of nonnegative real numbers such that there exists an integer $N_0 \geq 2$ for which

$$\varepsilon(N) \leq r_1 \, \varepsilon\left(\lfloor N/2 \rfloor\right) + r_2 \, \varepsilon\left(\lceil N/2 \rceil\right) + \delta(N) \qquad \text{for all } \; N \geq N_0, \tag{4.10}$$

where $r_1$, $r_2$ are nonnegative constants satisfying $r := r_1 + r_2 < 1$ and $\lim_N \delta(N) = 0$. Then,

$$\lim_N \varepsilon(N) = 0. \tag{4.11}$$

*Proof.* We have that $\limsup \varepsilon(N) = \limsup \varepsilon\left(\lfloor N/2 \rfloor\right) = \limsup \varepsilon\left(\lceil N/2 \rceil\right)$. Thus, by taking lim sup in both sides of (4.10) we get

$$0 \leq \limsup \varepsilon(N) \leq r_1 \limsup \varepsilon(N) + r_2 \limsup \varepsilon(N) + \limsup \delta(N),$$

i.e.

$$0 \leq \limsup \varepsilon(N) \leq r \limsup \varepsilon(N). \tag{4.12}$$

Since $r < 1$, to finish the proof it suffices to show that $\varepsilon(N)$ is bounded, and we can, e.g., see that as follows. Let $\delta^* := \max_N |\delta(N)|$ and

$$M := \max\left\{ \frac{\delta^*}{1 - r}, \; \max_{1 \leq N \leq N_0} |\varepsilon(N)| \right\}. \tag{4.13}$$

Then (noticing that (4.13) implies that $rM + \delta^* \leq M$), easy induction on $N$ implies that $\varepsilon(N) \leq M$ for all $N \geq 1$. ■

If $r_1 + r_2 = 1$, then (4.11) does not hold. For example, if $\varepsilon(N) = r_1 \, \varepsilon\left(\lfloor N/2 \rfloor\right) + (1 - r_1) \, \varepsilon\left(\lceil N/2 \rceil\right)$ for $N \geq 2$, then $\varepsilon(N) = \varepsilon(1)$ for all $N$.

**Corollary A1.** Let $z(N)$, $N = 1, 2, \ldots$, be a sequence of complex numbers satisfying the recursion

$$z(N) = \rho_1(N) \; z\left(\lfloor N/2 \rfloor\right) + \rho_2(N) \; z\left(\lceil N/2 \rceil\right) + b(N), \qquad N \geq 2, \tag{4.14}$$

where $\rho_1(N)$, $\rho_2(N)$, and $b(N)$ are complex sequences such that

$$\lim_N \rho_1(N) = \rho_1 \in \mathbb{C}, \quad \lim_N \rho_2(N) = \rho_2 \in \mathbb{C}, \quad \text{and} \quad \lim_N b(N) = b \in \mathbb{C}, \tag{4.15}$$

with $|\rho_1| + |\rho_2| < 1$. Then

$$\lim_N z(N) = \frac{b}{1 - \rho_1 - \rho_2}. \tag{4.16}$$

*Proof.* Pick $\epsilon > 0$ so that $|\rho_1| + |\rho_2| + 2\epsilon < 1$ and then choose $N_0 \geq 2$ so that $|\rho_1(N) - \rho_1| < \epsilon$ and $|\rho_2(N) - \rho_2| < \epsilon$ for all $N \geq N_0$. Then (4.14) implies

$$|z(N)| \leq (|\rho_1|+\epsilon) \left|z\left(\lfloor N/2 \rfloor\right)\right| + (|\rho_2|+\epsilon) \left|z\left(\lceil N/2 \rceil\right)\right| + b^*, \quad N \geq N_0, \tag{4.17}$$

where $b^* := \sup_N |b(N)| < \infty$, and the argument at the end of the proof of Lemma A1 applies to (4.17) and implies that the sequence $z(N)$ is bounded. Next, we write (4.14) as

$$z(N) - \frac{b}{1-\rho_1-\rho_2} = \rho_1 \left[z\left(\lfloor N/2 \rfloor\right) - \frac{b}{1-\rho_1-\rho_2}\right] + \rho_2 \left[z\left(\lceil N/2 \rceil\right) - \frac{b}{1-\rho_1-\rho_2}\right] + \delta(N), \tag{4.18}$$

where

$$\delta(N) := \left[\rho_1(N) - \rho_1\right] z\left(\lfloor N/2 \rfloor\right) + \left[\rho_2(N) - \rho_2\right] z\left(\lceil N/2 \rceil\right) + \left[b(N) - b\right]. \tag{4.19}$$

Notice that our assumptions for $\rho_1(N)$, $\rho_2(N)$, and $b(N)$, together with the fact that $z(N)$ is bounded, imply that $\lim_N \delta(N) = 0$. Thus, if we take absolute values in (4.18) and set

$$\varepsilon(N) := \left| z(N) - \frac{b}{1-\rho_1-\rho_2} \right|, \tag{4.20}$$

we obtain

$$\varepsilon(N) \leq |\rho_1|\, \varepsilon\left(\lfloor N/2 \rfloor\right) + |\rho_2|\, \varepsilon\left(\lceil N/2 \rceil\right) + |\delta(N)|, \qquad N \geq 2, \tag{4.21}$$

hence (since $|\rho_1| + |\rho_2| < 1$) Lemma A1 implies that $\lim_N \varepsilon(N) = 0$. ■

Corollary A1 is used in the proof of Theorem 6.
Finally, by using the above arguments we can also show the following simpler statement which is used in the proof of Theorem 3.

**Corollary A2.** Let $z_n$, $n = 0, 1, \ldots$, be a sequence of complex numbers satisfying the recursion

$$z_n = \rho z_{n-1} + b_n, \qquad n = 1, 2, \ldots, \tag{4.22}$$

where $|\rho| < 1$ and $\lim_n b_n \in \mathbb{C}$. Then,

$$\lim_n z_n = \frac{\lim_n b_n}{1-\rho}. \tag{4.23}$$

**Acknowledgments.** The author wishes to thank the two anonymous referees of another version of this article, as well as Professor Tom Britton, associate editor of JOMB, for their constructive comments and suggestions.